\documentclass{article}%
\usepackage{amsmath}
\usepackage{graphicx}
\usepackage{amsfonts}
\usepackage{amssymb}%
\providecommand{\U}[1]{\protect\rule{.1in}{.1in}}
\newtheorem{theorem}{Theorem}

\newtheorem{proposition}[theorem]{Proposition}

\begin{document}

\title{Some aspects of profinite group theory}
\author{Dan Segal}
\maketitle

\section{Introduction}

\subsection{Origins}

Profinite groups were introduced in number theory early in the last century.
First of all, the group of $p$-adic integers $\mathbb{Z}_{p}$ appeared as a
means for studying congruences: one can replace infinitely many congruences of
the form
\[
f(\mathbf{X})\equiv0\,(\operatorname{mod}p^{n})
\]
by a single \emph{equation}
\[
f(\mathbf{X})=0
\]
over $\mathbb{Z}_{p}$. There are two advantages to this approach. One is that
we can do arithmetic in a nice integral domain of characteristic zero, instead
of the messy finite rings $\mathbb{Z}/p^{n}\mathbb{Z}$. More importantly,
though, from a methodological point of view, what we have here is a technology
for replacing infinitely many hypotheses (about disparate small objects) with
a single hypothesis (about one large object): the ``large object'' -- the
$p$-adic integers in this case -- can then be studied by methods of algebra or
arithmetic. This process of ``mathematical reification'' is of course quite
traditional (as in the construction of the complex numbers), but is a
particularly characteristic feature of 20th century mathematics (Hilbert
space, representable functors,...).

As a profinite group, of course, $\mathbb{Z}_{p}$ is rather trivial, and its
main role in this context is as a \emph{ring}. Profinite groups of (much)
greater complexity were introduced by Krull. His insight was that the Galois
group of an infinite algebraic Galois extension of fields is in a natural way
a profinite group: it is a compact topological group, whose structure is
completely determined by the \emph{finite} Galois groups of all the finite
Galois subextensions. This led to the elegant modern formulations of class
field theory by Chevalley, Artin and Tate.

Later, Grothendieck introduced profinite groups into algebraic geometry, as
the fundamental groups of schemes. I shall say no more about these topics,
which are well beyond my competence: instead, I will concentrate on profinite
groups as objects of study for group theorists. This is not going to be a
comprehensive survey even of this limited subject: my intention is merely to
point the reader to some areas where interesting developments have recently
taken place, and that I happen to know something about (and the different
amounts of space devoted to the various topics in no way reflects their
relative importance).

I would like to acknowledge a special debt to both Alex Lubotzky and Avinoam
Mann, from whom -- way back in the last century -- I learnt new ways to do
group theory. Thanks are also due to Alex Lubotzky and Derek Holt for useful
contributions to this article.

\subsection{``Profinite group theory''}

This phrase, like `algebraic number theory', has a useful ambiguity (which is
lost on translation into French); I intend the reader to keep in mind both
meanings -- `the theory of profinite groups' and `a profinite approach to
group theory'.

A profinite group is what you get when you look at a (suitably coherent)
collection of finite groups all at once. In this context, `coherent' means
that the groups in question form an \emph{inverse system}: a family of finite
groups $(G_{\lambda})$ indexed by a directed set $\Lambda$, and for each pair
$\alpha,\,\beta\in\Lambda$ with $\alpha\leq\beta$ a homomorphism
$\theta_{\beta\alpha}:G_{\beta}\rightarrow G_{\alpha}$. Whenever $\alpha
\leq\beta\leq\gamma$ we require that $\theta_{\beta\alpha}\circ\theta
_{\gamma\beta}=\theta_{\gamma\alpha},$ and each $\theta_{\alpha\alpha}$ is the
identity automorphism. (To say that $\Lambda$ is a directed set means that
$\Lambda$ is partially ordered and that for every pair $\alpha,\,\beta
\in\Lambda$ there exists $\gamma\in\Lambda$ with $\gamma\geq\alpha$ and
$\gamma\geq\beta$.) The \emph{inverse limit} of this system, denoted
\[
\underset{\underset{\Lambda}{\longleftarrow}}{\lim}G_{\lambda},
\]
may be defined by a suitable universal property, or more concretely as a
subgroup $G$ of the Cartesian product of all the $G_{\lambda},$ as follows:
\begin{equation}
G=\left\{  (g_{\lambda})\mid\theta_{\beta\alpha}(g_{\beta})=g_{\alpha}\text{
whenever }\beta>\alpha\right\}  \leq\prod_{\lambda\in\Lambda}G_{\lambda}.
\label{invlim}%
\end{equation}
Thus $G$ maps naturally into each of the finite groups $G_{\lambda}$ (by
projecting to a factor), and $G$ is completely determined by the system
$\left(  G_{\lambda}\right)  _{\lambda\in\Lambda}$; the homomorphisms
$\theta_{\beta\alpha}$ are supposed to be included as part of the definition
of the system. So far, we have done little more than introduce a notation. The
key observation now is that $G$ is in a natural way a \emph{topological
group}: giving each of the finite groups $G_{\lambda}$ the discrete topology
we endow $\prod_{\lambda\in\Lambda}G_{\lambda}$ with the product topology;
instead of being discrete, this is a compact Hausdorff space, by Tychonoff's
Theorem. It is easy to see that the inverse limit $G$ is a closed subgroup, so
in this way $G$ becomes a \emph{compact Hausdorff topological group}. For each
$\lambda\in\Lambda$ the kernel $K_{\lambda}$ of the projection $\pi_{\lambda
}:G\rightarrow G_{\lambda}$ is an open normal subgroup of $G$, and the family
$\left\{  K_{\lambda}\right\}  $ forms a base for the neighbourhoods of $1$ in
$G$. In most naturally-arising situations, the maps $\theta_{\beta\alpha}$ are
all surjective, in which case one speaks of a \emph{surjective inverse
system}. The projections $\pi_{\lambda}$ are then also surjective ([RZ2],
Chapter 1), so the original groups $G_{\lambda}$ appear as continuous finite
images of $G$. To save repetition, I will make the blanket assumption that
\emph{all inverse systems under discussion are surjective}; this is no real
loss of generality, since to the inverse system $(G_{\lambda})$ we may
associate another one $(\pi_{\lambda}(G))$ which is surjective and has the
same inverse limit $G$.

The bare algebraic structure of $G$ may carry little information about the
original system of finite groups, but in combination with the topology it
closely reflects many properties of that system. Vaguely speaking, properties
of the topological group $G$ reflect \emph{uniform} properties of the groups
$G_{\lambda}$. For example, we find that $G$ is finitely generated (as a
\emph{topological} group) if and only if there exists $d\in\mathbb{N}$ such
that each of the groups $G_{\lambda}$ can be generated by $d$ elements; more
subtle relationships of this kind will be discussed below.

Certain classes of profinite groups have special names: if the finite groups
$G_{\lambda}$ all belong to some class of groups $\mathcal{C}$ (and the
inverse system is surjective), then $G=\underset{\longleftarrow}{\lim
}G_{\lambda}$ is called a pro-$\mathcal{C}$ group. When $\mathcal{C}$ is the
class of finite $p$-groups for some prime $p$ one calls $G$ a pro-$p$ group.

In practice, most of the questions studied in ``profinite group theory'' arise
in one of the following contexts, which are not mutually exclusive.

(1) Questions about some naturally-defined family of finite groups, for
example finite $p$-groups; see \S \ref{p-groups}.

(2) Questions about infinite groups that can be approached through their
profinite completions; this may be construed as a subcase of (1) where the
family of finite groups consists of finite quotients of some fixed infinite
group. See \S \S \ref{local}, \ref{sg}, \ref{csp}.

(3) Questions about profinite groups as such; these may be analogues in the
profinite category of familiar group-theoretic questions (\S \ref{free}%
,\thinspace\ref{fp}), they may arise from number theory and field theory via
the Galois group (see [FJ], [B], [dSF]), or they may be a new kind of question
specific to the profinite situation (\S \S \ref{fg}, \ref{prob}). As we shall
see, such investigations sometimes lead to new results about abstract groups,
finite or infinite.

For definitions and the basic properties of profinite groups, consult [DDMS]
Chapter 1, [W2] Chapter 1, or [RZ2] Chapter 2. Each of these books goes on to
discuss various specific topics; some of them are mentioned below, but all are
worth studying. Galois-theoretic applications of profinite groups are pursued
at length in [FJ]. Various aspects of pro-$p$ groups are discussed in detail
in [NH], which includes a substantial list of open problems.

The first substantial treatment in book form of profinite groups was Serre's
influential book [CG]: as the title suggests, this is primarily concerned with
homological matters and is slanted towards number theory.

\section{Local and global\label{local}}

An important strand in number theory is the investigation of so-called
`local-global' principles. A typical question is like this: suppose a certain
Diophantine equation $f(\mathbf{X})=c$ can be solved modulo $m$ for every
$m\in\mathbb{N},$ does it follow that $f(\mathbf{X})=c$ has a solution in
integers? The most famous example is the Hasse-Minkowski Theorem, which gives
an affirmative answer for \emph{rational} solutions, at least, when
$f(\mathbf{X})$ is an indefinite quadratic form. A fruitful way of formulating
such statements was introduced by Hasse (inspired by Hensel): instead of
considering congruences one takes the (equivalent) hypothesis that
$f(\mathbf{X})=c$ is solvable in every $p$-adic field (if we call $\mathbb{R}$
the `$\infty$-adic field' this also subsumes the condition that the quadratic
form be indefinite). In this case, the equation is said to be solvable
`locally'; if this implies the existence of a rational solution one has a
`local-global' theorem, and the equation is said to satisfy the `Hasse Principle'.

In the case of quadratic forms, the answer for \emph{integral} solutions is a
little more complicated: it may be that $f(\mathbf{X})=c$ is not solvable in
integers, but at least we can say that $f^{\prime}(\mathbf{X})=c$ is solvable
where $f^{\prime}(\mathbf{X})$ is one of \emph{finitely many} quadratic forms,
that constitute the \emph{genus} of $f$.

What has this to do with group theory? Thinking of the above as the search for
properties of $\mathbb{Z}$ that are determined by properties of the collection
of finite rings $\mathbb{Z}/m\mathbb{Z},$ we can generalize as follows: to
what extent are properties of an infinite group $G$ determined by the finite
quotient groups of $G$? This is a natural enough question in itself; it also
has a further philosophical motivation, connected with \emph{decision
problems}. The point is that a `local-global' theorem in group theory, as in
number theory, often implies a corresponding decidability theorem. Rather than
stating this as a formal metatheorem let me illustrate with an example which
should make the idea clear. A group $G$ is \emph{conjugacy separable} if it
has the following property: if $x,\,y\in G$ are such that $\pi(x)$ is
conjugate to $\pi(y)$ in $\pi(G)$ for every homomorphism $\pi$ from $G$ to any
finite group, then $x$ and $y$ are conjugate in $G$; this is the
`local-global' property of conjugacy in $G$. On the other hand, $G$ has
\emph{solvable conjugacy problem} if there is a uniform algorithm that
decides, given any two elements of $G$, whether or not they are conjugate in
$G$.

\begin{theorem}
Every finitely presented conjugacy separable group has solvable conjugacy problem.
\end{theorem}

\noindent The algorithm consists of two procedures, run simultaneously. The
first one lists all consequences of the relations in a given presentation of
$G$, while the second one enumerates all homomorphisms $\pi$ from $G$ to
finite groups, and for each such $\pi$ lists the (finitely many) pairs of
non-conjugate elements in $\pi(G)$. Now given $x$ and $y\in G$, we run both
procedures until \emph{either} the first one spits out an equality $x^{g}=y$
\emph{or} the second one spits out a pair $(\pi(x),\pi(y))$. In the first case
we conclude that $x$ and $y$ are conjugate, in the second that they are not;
the hypothesis that $G$ is conjugacy separable ensures that one or other of
the cases must arise.

Of course, no sane person would try to implement such a stupid algorithm; its
interest is theoretical. It shows that combinatorial group theorists shouldn't
waste their time trying to prove the unsolvability of the conjugacy problem in
the case of conjugacy separable groups. The same applies to the \emph{word}
problem in \emph{residually finite} groups: a group is residually finite if
its subgroups of finite index intersect in $\{1\}$, which is equivalent to
saying that any two distinct elements have distinct images in at least one
finite quotient of the group -- the local-global property for equality of
elements (while the `word problem' asks for an algorithm to determine equality
of group elements, given as words on a fixed generating set).

As in number theory, there is a useful reformulation for group-theoretic
`local-global' questions. The family of all finite quotients of $G$ naturally
forms an inverse system of finite groups, with respect to the quotient maps
\[
G/N\rightarrow G/M
\]
where $N\leq M$ are normal subgroups of finite index in $G$. The inverse limit
of this system is the \emph{profinite completion }$\widehat{G}$ of $G$; and
the question becomes: what properties of $G$ are determined by properties of
the profinite group $\widehat{G}$?

The family of quotient maps $G\rightarrow G/N$ induces a natural homomorphism
$\iota:G\rightarrow\widehat{G}$. The kernel of $\iota$ is the \emph{finite
residual} $R(G)$ of $G$, which is the intersection of all subgroups of finite
index in $G$. Evidently $\widehat{G}=\widehat{G/R(G)}$, so knowledge of
$\widehat{G}$ will at best give us information about $G/R(G)$; thus it is
sensible to restrict attention to groups $G$ for which $R(G)=1$, that is
residually finite groups. If $G$ is residually finite then the map $\iota$ is
injective, and we use it to identify $G$ with a subgroup of $\widehat{G}$.
This amounts to identifying an element $g\in G$ with the `diagonal' element
\[
(gN)_{N\in\mathcal{N}}\in\widehat{G}\leq\prod_{N\in\mathcal{N}}G/N
\]
where $\mathcal{N}$ is the family of all normal subgroups of finite index in
$G$.

To say that $G$ is residually finite, then, amounts to saying that two
elements of $G$ are equal if and only if they map to equal elements of
$\widehat{G}$. Similarly, $G$ is conjugacy separable if and only if for pairs
of elements of $G$, conjugacy in $\widehat{G}$ implies conjugacy in $G$; this
is equivalent to saying that \emph{each conjugacy class} in $G$ \emph{is
closed} in the profinite topology of $G$, that is, the topology induced from
$\widehat{G}$, in which a base for the neighbourhoods of $1$ in $G$ is given
by the family $\mathcal{N}$ of all normal subgroups of finite index in $G$.
(Analogously, $G$ is residually finite if and only if \emph{points} are closed
in $G$ -- in this context this is equivalent to the Hausdorff property )

Well known classes of residually finite groups include the \emph{free groups}
and the \emph{virtually polycyclic groups}. (A group is \emph{virtually} $P$
if it has a normal $P$-subgroup of finite index.) In fact, groups in these
classes have many good local-global properties: in particular, they are

\begin{itemize}
\item \emph{conjugacy separable} and

\item \emph{subgroup separable};
\end{itemize}

\noindent a group $G$ is subgroup separable if every finitely generated
subgroup is closed in the profinite topology of $G$; this is equivalent to
saying that for each finitely generated subgroup $H,$ the property `being in
$H$' is a local-global property of elements of $G$. This has the important
consequence that the closure of $H$ in $\widehat{G}$ is naturally isomorphic
to $\widehat{H}$.

That free groups are subgroup separable was proved by Marshall Hall in 1949
(see [LS], Chapter 1, Prop. 3.10). The fact that virtually polycyclic groups
are subgroup separable is quite elementary (see [S], Chapter 1). The fact that
they are conjugacy separable, however, depends on an interesting result in
algebraic number theory (due to F. K. Schmidt and Chevalley):

\begin{theorem}
\label{chev}Let $\mathcal{O}$ be the ring of integers in an algebraic number
field, with group of units $\mathcal{O}^{\ast}$. Then every subgroup of finite
index in $\mathcal{O}^{\ast}$ contains a `congruence subgroup'
$(1+m\mathcal{O})\cap\mathcal{O}^{\ast}$ $\ (0\neq m\in\mathbb{Z})$.
\end{theorem}

\noindent One says that $\mathcal{O}^{\ast}$ has the \emph{congruence subgroup
property}: equivalently, the profinite topology on the \emph{additive} group
of $\mathcal{O}$ induces the profinite topology on the \emph{multiplicative}
group $\mathcal{O}^{\ast}$. More generally, let $M$ be any virtually
polycyclic group and $G$ a group of automorphisms of $M$. We may define a
topology on $G$ by choosing as a base for the neighbourhoods of $1$ the family
of subgroups
\[
\mathrm{C}_{G}(M/M^{m})\,\,(m\in\mathbb{N});
\]
this is the \emph{congruence topology} on $G$.

\begin{theorem}
\label{congtop}Let $M$ be a virtually polycyclic group and $G$ a virtually
polycyclic subgroup of $\mathrm{Aut}(M)$. Then

\emph{(i)} Each orbit $a^{G}$ ($a\in M$) is closed in the profinite topology
of $M$;

\emph{(ii)} $G$ is closed in the congruence topology of $\mathrm{Aut}(M);$

\emph{(iii)} the congruence topology on $G$ is the same as the profinite topology.
\end{theorem}

Theorem \ref{chev} is the special case of (iii) where $M=\mathcal{O}$ and
$G=\mathcal{O}^{\ast}$ (acting by multiplication). Part (iii) follows directly
from (ii) applied to arbitrary subgroups of finite index in $G$. Part (ii)
follows from (i) applied to the orbit of $(a_{1},\ldots,a_{d})$ in the group
$M^{(d)}$, where $\left\{  a_{1},\ldots,a_{d}\right\}  $ is a generating set
for $M$ and $G$ acts diagonally. Part (i), a generalized version of conjugacy
separability, may be reduced to an application of Theorem \ref{chev} by
`d\'{e}vissage', arguing by induction on the Hirsch length of $M$. See [S],
Chapter 4 (I assumed there that $M$ is free abelian, but the general case is
no harder). A further generalization is given in [S1], \S 8.

There is an essentially equivalent formulation of (i) in terms of
\emph{derivations}: a derivation from $G$ to $M$ is a map $\delta:G\rightarrow
M$ such that $\delta(xy)=\delta(x)^{y}.\delta(y)$ for all $x,\,y\in G$ (such
maps are also called \emph{crossed homomorphisms} or $1$\emph{-cocyles}).
Among these are the \emph{inner derivations} $\delta_{a}:x\mapsto a^{x}a^{-1}$
($a\in M$ fixed). Since $a^{G}=\delta_{a}(G)\cdot a$ we see that (i) is a
special case of

\begin{theorem}
\label{closedder}Let $M$ and $G$ be virtually polycyclic groups, with $G$
acting on $M$. If $\delta:G\rightarrow M$ is a derivation then the set
$\delta(G)$ is closed in the profinite topology of $M$.
\end{theorem}

The action of $G$ on $M$ induces an action of $\widehat{G}$ on $\widehat{M}$,
and a derivation $\delta:G\rightarrow M$ induces a continuous derivation
$\widehat{\delta}:\widehat{G}\rightarrow\widehat{M}$. One may deduce

\begin{theorem}
\label{H1}Let $G$ and $M$ be as above. Then the natural mapping
\[
H^{1}(G,M)\rightarrow H^{1}(\widehat{G},\widehat{M})
\]
is injective.
\end{theorem}

\noindent Here, $H^{1}(G,M)$ is the `non-abelian cohomology' set defined in
[CG], Chapter 1. Another application of Theorem \ref{closedder} gives

\begin{proposition}
\label{der}Let $a\in M$. Then $\widehat{\delta}^{-1}(a)$ is equal to the
closure of $\delta^{-1}(a)$ in $\widehat{G}$.
\end{proposition}

\noindent This applies in particular when $\delta$ is a homomorphism, and
shows that the functor $G\mapsto\widehat{G}$ is \emph{exact }on virtually
polycyclic groups. This can also be seen by a direct elementary argument; but
the following excellent properties of this functor depend on the full strength
of Proposition \ref{der}:

\begin{theorem}
\label{blob}Let $G$ be a virtually polycyclic group and $H,\,K$ subgroups of
$G$. Then
\begin{align*}
\mathrm{C}_{\overline{K}}(\overline{H})  &  =\overline{\mathrm{C}_{K}(H)}\\
\mathrm{N}_{\overline{K}}(\overline{H})  &  =\overline{\mathrm{N}_{K}(H)}\\
\overline{H}\cap\overline{K}  &  =\overline{H\cap K}.
\end{align*}
where $\overline{X}$ denotes the closure of a set $X$ in $\widehat{G}$.
\end{theorem}

\noindent See [RSZ], \S 2. This is applied together with the geometric study
of profinite groups acting on `profinite trees' to establish

\begin{theorem}
\label{csep}\emph{[RSZ] }Let $G$ a group that is obtained from virtually free
groups and virtually polycyclic groups by forming finitely many successive
free products, amalgamating cyclic subgroups. Then $G$ is conjugacy separable.
\end{theorem}

Sometimes, a purely global result can be deduced by a `local-global' argument
from the finite case: this is how Hirsch (in 1954) proved that the Frattini
subgroup of a polycyclic group is nilpotent. As an application of Theorem
\ref{blob}, we show in [NS3] that if $G$ is a virtually polycyclic group, $N$
is a normal subgroup, and $G$ is isomorphic to $N\times G/N$, then $N$ is
actually a direct factor of $G$; this is deduced from the special case of
finite groups, a recent theorem due to J. Ayoub.

Another famous decision problem in group theory is the \emph{isomorphism
problem}: to decide, given two finite group presentations, whether or not they
define isomorphic groups. Suppose $\mathcal{C}$ is a class of groups having
the `local-global property for isomorphism' -- that is, for $G$ and $H$ in
$\mathcal{C}$ one has $G\cong H$ if and only if $\mathcal{F}(G)=\mathcal{F}%
(H)$, where $\mathcal{F}(G)$ denotes the set of isomorphism types of finite
quotient groups of $G$. Then it is easy to see, by a modification of the
argument above, that the isomorphism problem for finitely presented groups in
$\mathcal{C}$ has a positive solution. Examples of such classes $\mathcal{C}$
are the finitely generated free groups and the finitely generated abelian
groups. As polycyclic groups are not so very different from finitely generated
abelian groups, one might wonder whether they, also, have the local-global
property for isomorphism. The answer is `no': examples demonstrating this are
given in [S], Chapter 11. Some of these examples are constructed using
integral quadratic forms that are `locally equivalent' but not equivalent over
$\mathbb{Z}$. Such quadratic forms, however, do belong to the same genus,
which consists of \emph{finitely many} integral equivalence classes. And this
finiteness property does indeed generalize:

\begin{theorem}
\label{gps}\emph{[GPS] }Given any set $\mathcal{X}$ of isomorphism types of
finite groups, there are at most finitely many isomorphism types of virtually
polycyclic groups $G$ such that $\mathcal{F}(G)=\mathcal{X}$.
\end{theorem}

\noindent The proof does not tell us exactly \emph{how many} isomorphism
types, so the theorem does \emph{not} imply a positive solution for the
isomorphism problem in this case. That requires other methods, and may be
found in [S1].\ While the statement of the theorem does not explicitly mention
profinite groups, it is clear (if $G$ is a finitely generated group!) that the
set $\mathcal{F}(G)$ both determines and is determined by (the
\emph{topological }group) $\widehat{G}$, so the result amounts to saying that
for virtually polycyclic groups, the profinite completion `determines the
group up to finitely many possibilities'. (In fact, in this case the
topological group $\widehat{G}$ is uniquely determined by its underlying
abstract group: see \S \ref{fg} below.)

Given any subgroup $G$ of a profinite group $P,$ the inclusion $G\rightarrow
P$ induces a natural continuous homomorphism $\pi:\widehat{G}\rightarrow P$.
This is surjective if and only if $G$ is dense in $P$; it is injective if and
only if the topology induced on $G$ as a subspace of $P$ is the profinite
topology of $G$, in which case we say that $G$ has the \emph{congruence
subgroup property}, or CSP\emph{\ }(by analogy with the special case
$P=\mathrm{Aut}(M)$ discussed above). Now we can reformulate Theorem \ref{gps} as

\begin{theorem}
\label{gpsprof}Let $P$ be a profinite group. Let $\mathcal{S}$ denote the set
of all virtually polycyclic subgroups that are dense in $P$ and have CSP. Then
$\mathcal{S}$ consists of finitely many orbits of $\mathrm{Aut}(P)$.
\end{theorem}

\noindent In fact, using Theorem \ref{ur}, below, one can show that if $P$ is
the profinite completion of a virtually polycyclic group, then any finitely
generated residually finite group $G$ with $\widehat{G}\cong P$ is itself
virtually polycyclic (by considering the dimension of the Sylow pro-$p$
subgroups of $\widehat{G}$: see \S \ref{pag}). So in Theorem \ref{gpsprof} we
can replace `virtually polycyclic' by `finitely generated' as long as we add
the hypothesis that $P$ contain at least one dense virtually polycyclic
subgroup with CSP.

The advantage (indeed, the necessity) of this `profinite' approach is apparent
as soon as one embarks on the proof of this theorem. One of the first steps,
for example, is to show that for $G\in\mathcal{S}$ the closure in $P$ of the
Fitting subgroup of $G$ is precisely the Fitting subgroup of $P$. Since $G$ is
subgroup separable it follows that $\widehat{\mathrm{Fit}(G)}$ is determined
by $\widehat{G}$; the problem can now be broken into two cases: (1) the case
of \emph{nilpotent} groups, (2) the study of groups $G$ for which not only
$\widehat{G}$ but also $\mathrm{Fit}(G)$ are fixed. Both parts are difficult,
and depend on deep results in the arithmetic theory of algebraic groups,
results that generalize classical finiteness properties of quadratic forms.
The key fact is the following analogue of Theorem \ref{congtop}(i):

\begin{theorem}
\emph{(Borel and Serre)} Let $\Gamma$ be an arithmetic group acting rationally
on $M=\mathbb{Z}^{d}$. Then each `local orbit' of $\Gamma$ in $M$ is the union
of finitely many orbits of $\Gamma$.
\end{theorem}

\noindent A \emph{local orbit} here means a set of the form $M\cap
a^{\widetilde{\Gamma}}$ where $\widetilde{\Gamma}$ is the \emph{integral adele
group} associated to $\Gamma,$ acting on $\widehat{M}$. For all this, see [S],
Chapters 9 and 10.

\bigskip

A much simpler question than that of isomorphism is the following: what is the
minimal size of a generating set for a group $G$ ? This number is denoted
$\mathrm{d}(G)$.

\begin{theorem}
\emph{[LW] }If $G$ is a virtually polycyclic group then $\mathrm{d}%
(G)\leq\mathrm{d}(\widehat{G})+1$.
\end{theorem}

\noindent Of course, $\mathrm{d}(\widehat{G})$ here denotes the minimal size
of a \emph{topological} generating set for $\widehat{G}$, so what the result
is saying is that if every finite quotient of $G$ can be generated by $d$
elements, then $G$ itself can be generated by $d+1$ elements. This is a hard
theorem due to Linnell and Warhurst. It is very easy to find cases where
$\mathrm{d}(G)=\mathrm{d}(\widehat{G})$ (abelian groups for example), and not
much harder to find cases where $\mathrm{d}(G)=\mathrm{d}(\widehat{G})+1$
(using a ring of algebraic integers that is not a PID). If we had an algorithm
for deciding whether a polycyclic group is of the first or of the second type,
we could then effectively determine $\mathrm{d}(G)$ for such groups $G$, by a
version of the `stupid double-enumeration procedure' described above. But --
as far as I know -- we don't. Indeed the following challenging problem is
still open (even for the `easy' case of virtually abelian groups!):

\medskip

\noindent\textbf{Problem.} \emph{Find an algorithm that determines}
$\mathrm{d}(G)$ \emph{for every polycyclic group} $G$.

\medskip

\noindent For the currently known decision procedures for polycyclic groups
see [BCRS], [S1], [E1] and [E2].

\bigskip

A uniform bound for $\mathrm{d}(\overline{G})$ over all the finite images
$\overline{G}$ of a group $G$ is just one example of what I call an `upper
finiteness condition': a uniform bound for some measure of size, or growth, on
all the finite quotients of a group. Any such condition certainly means
something for the global structure of a group, and the challenge is to find
out what it is. This programme is discussed in Section \ref{sg} below.

\section{$p$-Adic analytic groups\label{pag}}

The theory of Lie groups is without doubt one of the central pillars of
twentieth-century mathematics (not to mention physics!). Quite early in the
century, an analogous theory of `$p$-adic Lie groups' received some attention:
these $p$\emph{-adic analytic groups }have the underlying structure of an
analytic manifold over the field $\mathbb{Q}_{p}$, and the group operations
are given locally by convergent $p$-adic power series. The global structure
and cohomology theory of $p$-adic analytic groups were elucidated by Michel
Lazard in a magisterial paper [L], published in 1965. One of his key
discoveries was that each compact $p$-adic analytic group has an open subgroup
(necessarily of finite index) which is a finitely generated pro-$p$ group, and
any pro-$p$ group arising in this situation has a certain special algebraic
property; conversely, every finite extension of a finitely generated pro-$p$
group with this property has, in a natural way, the structure of a compact
$p$-adic analytic group.

The `special property' discovered by Lazard is that of being \emph{powerful},
a term introduced later by Lubotzky and Mann in [LM]. The pro-$p$ group $G$ is
powerful if $G/G^{p}$ is abelian (when $p=2$ we require that $G/G^{4}$ be
abelian). Thus powerful groups are `abelian to a first approximation', and
Lubotzky and Mann went on to show that in fact such groups resemble abelian
groups in many ways: for example, in a $d$-generator powerful group every
closed subgroup can be generated by $d$ elements. Thus such a group has
\emph{finite rank}, where the rank of a profinite group $G$ is defined by
\[
\mathrm{rk}(G)=\sup\{\mathrm{d}(H)\mid H\leq_{c}G\}
\]
(here $\mathrm{d}(H)$ denotes the minimal size of a (topological) generating
set for $H$, and $H\leq_{c}G$ means `$H$ is a closed subgroup of $G$').
Conversely, they proved that every pro-$p$ group of finite rank has an open
(hence of finite index) powerful subgroup. With Lazard's result, this shows
that \emph{a pro-}$p$\emph{\ group is }$p$\emph{-adic analytic if and only if
it has finite rank}.

This opened the way to a more group-theoretic approach to the whole topic,
expounded in detail in the book [DDMS]. The resulting theory has found
numerous applications. Applications to finite $p$-groups are discussed in
\S \ref{p-groups} below. Many applications to infinite group theory are based
on Lubotzky's observation that a compact $p$-adic analytic group is a linear
group over $\mathbb{Q}_{p}$, by Ado's Theorem: it should be mentioned that the
correspondence Lie groups $\leftrightarrow$ Lie algebras works even better in
the $p$-adic case than in the classical case. This leads to the `Lubotzky
linearity criterion', see \S \ref{fg}. It implies that any infinite group
which is residually a finite $p$-group and whose pro-$p$ completion has finite
rank is in fact a linear group; such a group can then be attacked with various
tools from linear group theory. A strikingly successful example of this
strategy is discussed in the following section.

Other group-theoretic applications are described in [DDMS]. Our relatively
good understanding of pro-$p$ groups of finite rank has encouraged the
investigation of wider classes of pro-$p$ groups, and this is currently a
lively area of research. Many recent developments are described in the book [NH].

\section{Upper finiteness conditions and \newline subgroup growth\label{sg}}

\subsection{`Upper finiteness conditions'}

Let us consider the implications for a group of imposing various restrictions
on its finite quotients.\medskip

\textbf{1. } The \emph{rank} $\mathrm{rk}(Q)$ of a finite group $Q$ is the
least integer $r$ such that every subgroup of $Q$ can be generated by $r$
elements. The \emph{upper rank} of any group $G$ is
\[
\mathrm{ur}(G)=\sup\left\{  \mathrm{rk}(Q)\mid Q\in\mathcal{F}(G)\right\}  .
\]
This is none other than the rank of $\widehat{G}$, defined above.

\begin{theorem}
\label{ur}\emph{[MS1]} Let $G$ be a finitely generated residually finite
group. Then $G$ has finite upper rank if and only if $G$ is virtually soluble
of finite rank.
\end{theorem}

An infinite group $G$ is said to have finite rank if there exists an integer
$r$ such that every \emph{finitely generated} subgroup of $G$ can be generated
by $r$ elements. Soluble groups of finite rank are quite easy to describe:
such a group that is also finitely generated and residually finite is a finite
extension of a triangular matrix group over a ring of the form $\mathbb{Z}%
[1/m]$. This theorem is making two remarkable assertions: (a) that `(bounded)
finite rank' is a local-global property, and (b) that a \emph{numerical bound}
(on the size of generating sets, in this case) implies a \emph{structural}
\emph{algebraic} property, namely solubility.\medskip

\textbf{2.} For a finite group $Q$, the \emph{number of subgroups} of $Q$ is
denoted $s(Q)$. A group $G$ has \emph{weak polynomial subgroup growth}, or
wPSG, if there exists a constant $\alpha$ such that
\begin{equation}
s(Q)\leq\left|  Q\right|  ^{\alpha} \label{alpha}%
\end{equation}
for every $Q\in\mathcal{F}(G)$.

\begin{theorem}
\emph{[LMS], [S2]} Let $G$ be a finitely generated residually finite group.
Then $G$ has wPSG if and only if $G$ is virtually soluble of finite rank.
\end{theorem}

\noindent The alert reader will have noticed that this theorem implies the
preceding one, since if $G$ has finite upper rank we can take $\alpha
=\mathrm{ur}(G)$ and deduce that $G$ has wPSG.\medskip

\textbf{3.} A group $G$ has \emph{polynomial index growth}, or PIG, if there
exists a constant $\alpha$ such that
\[
\left\vert Q\right\vert \leq(\exp Q)^{\alpha}%
\]
for every $Q\in\mathcal{F}(G)$, where $\exp Q$ denotes the exponent of $Q$.
This is equivalent to saying that $\left\vert Q/Q^{m}\right\vert \leq
m^{\alpha}$ for every $Q\in\mathcal{F}(G)$ and every $m\in\mathbb{N}$. It is
easy to see that every soluble group of finite rank has PIG, but the converse
is far from true: Balog, Mann and Pyber [BMP] construct a finitely generated
residually finite group with PIG which has finite simple quotients of
unbounded ranks. However, if we \emph{assume} solubility we have

\begin{theorem}
\emph{[PS]} Let $G$ be a finitely generated soluble residually finite group.
Then $G$ has PIG if and only if $G$ has finite rank.
\end{theorem}

PIG and other upper finiteness conditions are discussed in detail in Chapter
12 of [SG] (where this last result appears as an open problem).

\subsection{Subgroup growth}

A group $G$ has `weak PSG' if it doesn't have very many subgroups of each
finite index. More generally, it is interesting to study just how many
subgroups there are of each index: that is, to study the function $n\mapsto
a_{n}(G)$ where $a_{n}(G)$ denotes the number of subgroups of index $n$ in $G
$. This function is well defined as long as $G$ is finitely generated. When
$G$ is a profinite group, $a_{n}(G)$ denotes the number of \emph{open}
subgroups of index $n$ in $G$, and again is well defined if $G$ is
(topologically) finitely generated. Moreover, it is easy to see that if $G$ is
any abstract group, then $a_{n}(G)=a_{n}(\widehat{G})$; in this sense,
\emph{subgroup growth} -- i.e. the behaviour of the function $n\mapsto
a_{n}(G)$ -- is a `profinite' property of groups.

A comprehensive account of this topic is given in the book [SG], where the
advantages of the `profinite philosophy' are amply illustrated; let me just
mention a few of the highlights, under three headings. We will denote by
$s_{n}(G)$ the number of subgroups (or open subgroups) of index \emph{at most}
$n$ in the group $G$.\medskip

\textbf{`Analytic problems': }what does a given restriction on the subgroup
growth imply for the algebraic structure of a group?

A group $G$ has \emph{polynomial subgroup growth}, or PSG, if $\log
s_{n}(G)=O(\log n)$. This obviously implies wPSG, and it is a deep result
(depending on CFSG) that the two conditions are in fact equivalent. Thus the
theorem stated above is equivalent to

\begin{theorem}
\label{psgthm}\emph{[LMS] }Let $G$ be a finitely generated residually finite
group. Then $G$ has PSG if and only if $G$ is virtually soluble of finite rank.
\end{theorem}

\noindent The difficult part is `only if'. The original proof of this (though
not the one presented in [SG]) starts by considering the pro-$p$ completions
of $G.$ Lubotzky and Mann proved that every pro-$p$ group with PSG is $p$-adic
analytic, from which it follows that if $G$ has PSG then $\widehat{G}_{p}$ is
a $p$-adic analytic group, and therefore \emph{linear}. Thus if $G$ happens to
embed into $\widehat{G}_{p}$ then $G$ itself is a linear group. One can then
use `Strong Approximation' results (specifically, Theorem \ref{LA} stated in
\S \ref{fg}, below) to reduce to the case of arithmetic groups, and the proof
is concluded by an explicit counting of congruence subgroups in such groups
(see \S \ref{csp}). In the general case, further arguments are required,
depending among other things on CFSG.

While a finitely generated residually finite group with PSG must be virtually
soluble, this is not true for finitely generated \emph{profinite} groups with
PSG. These are characterized in [SSh]: such a profinite group is (virtually)
an extension of a prosoluble group of finite rank by the Cartesian product of
a family of finite quasisimple groups of Lie type satisfying certain very
precise arithmetical conditions. (In view of the preceding theorem, such a
group can only be the profinite completion of a finitely generated abstract
group in the special case where this family of quasisimple groups is
\emph{finite}.)

Like much of `pure' profinite group theory, the characterization of profinite
groups with PSG quickly reduces to a problem of \emph{finite} group theory:
establishing uniform bounds for several structural parameters of a finite
group $G$ in terms of the parameter $\alpha$ defined in (\ref{alpha}), above.
The same applies to many other results that relate the algebraic structure of
a profinite group to its rate of subgroup growth, when this is faster than polynomial.

.\medskip

\textbf{`Synthetic problems':} under this heading comes the problem of
constructing groups that demonstrate particular types of subgroup growth. A
group $G$ is said to have \emph{growth type} $f$ if
\begin{align*}
\log s_{n}(G)  &  =O(\log f(n))\\
\log s_{n}(G)  &  \neq o(\log f(n)).
\end{align*}
It is not difficult to construct finitely generated profinite groups with
more-or-less arbitrary growth type, by forming Cartesian products of suitable
collections of finite groups [MS2]. To do the same for finitely generated
\emph{abstract} groups is much harder, but we have

\begin{theorem}
\emph{([P], [S4])} Let $g:\mathbb{N\rightarrow R}_{+}$ be a `good'
non-decreasing function with $g(n)=O(n)$. Then there exists a finitely
generated group $G$ having growth type $n^{g(n)}$.
\end{theorem}

\noindent The condition `good' here is a mild restriction of a technical
nature, that need not concern us. The bound $g(n)=O(n)$ is necessary, because
the fastest possible growth type for any finitely generated group is easily
seen to be $n^{n}$. Thus the point of the theorem is that essentially every
`not impossible' growth type is actually exhibited by some finitely generated group.

The proof is in two stages. The first is to construct a suitable profinite
group $P$ with the specified growth type; the second, harder part, is to show
that this $P$ is the profinite completion of some finitely generated abstract
group (this is what `suitable' means here: the easy groups given in [MS2]
don't have this property). That is, we require $P$ to contain a dense finitely
generated subgroup $G$ that has the \emph{congruence subgroup property}, as
defined in \S \ref{local}, above. In fact two different constructions are
used: when $g(n)=O(\log\log n)$ one takes $P$ to be a certain group of
automorphisms of a rooted tree; this construction is discussed in \S \ref{fg},
below. When $\log n=O(g(n))$ one takes $P$ to be the Cartesian product of a
suitable family of finite alternating groups; in this case, the dense subgroup
$G$ does not quite have the CSP, but close enough: it turns out that the
kernel of the natural epimorphism $\widehat{G}\rightarrow P$ is a procyclic
group, which is enough to ensure that $G$ has the same subgroup growth type as
$P$. For full details see Chapter 13 of [SG] (a different and more general
construction has recently been obtained in [KN]; see Theorem \ref{KN} below).

\medskip

`\textbf{Zeta functions': \ }Having associated to a finitely generated group
$G$ the numerical sequence $(a_{n}(G))$, it is natural to wonder about the
arithmetical properties of this sequence. The `growth type' defined above is
one crude measure, but can we obtain more refined information? This question
has been studied in depth for certain types of groups: (a) free groups,
one-relator groups and free products of finite groups, (b) finitely generated
nilpotent groups, and (c) $p$-adic analytic pro-$p$ groups.

I will say no more about the class (a). This is the subject of many papers by
Thomas M\"{u}ller, using methods of combinatorics and analysis; for references
and some sample results see Chapter 14 of [SG]. Groups of types (b) and (c)
have polynomial subgroup growth: in this case, it is convenient to encode the
sequence $a_{n}(G)$ in a generating function
\[
\zeta_{G}(s)=\sum_{n=1}^{\infty}a_{n}(G)n^{-s}
\]
where $s$ is a complex variable. This `zeta function' represents a complex
analytic function, regular on some half-plane $\operatorname{Re}(s)>\alpha$;
here the abscissa of convergence $\alpha$ is given by
\[
\alpha=\inf\left\{  \gamma\mid s_{n}(G)=O(n^{\gamma})\right\}  ,
\]
a finite number when $G$ has PSG.

For a fixed nilpotent group $G$, it is easy to see that the arithmetical
function $a_{n}(G)$ is multiplicative, i.e. if $m$ and $n$ are coprime then
$a_{mn}(G)=a_{m}(G)a_{n}(G)$. This implies the `Euler product' decomposition
\[
\zeta_{G}(s)=\prod_{p}\zeta_{G,p}(s)
\]
where the product is over all primes and the `local factors' are defined by
\[
\zeta_{G,p}(s)=\sum_{j=0}^{\infty}a_{p^{j}}(G)p^{-js}.
\]
We showed in [GSS] that when $G$ is a finitely generated nilpotent group, for
each prime $p$ the series $\zeta_{G,p}(s)$ represents a\emph{\ rational
function }in $p^{-s}$ (with rational coefficients); to see why this is
reasonable, note that when $G$ is the infinite cyclic group $\zeta_{G}$ is the
Riemann zeta function, and $\zeta_{G,p}(s)=\frac{1}{1-p^{-s}}$. The proof
applies a general theorem about $p$-adic integrals, proved by Denef using
methods of $p$-adic model theory. Now, still assuming that $G$ is finitely
generated and nilpotent, we have $\zeta_{G,p}(s)=\zeta_{P}(s)$ where
$P=\widehat{G}_{p}$ is the pro-$p$ completion of $G$; and $P$ in this case is
a $p$-adic analytic pro-$p$ group. Thus the rationality theorem just mentioned
is a very special case of

\begin{theorem}
\emph{[dS1] }If $P$ is a compact $p$-adic analytic group then $\zeta_{G,p}(s)$
is a rational function over $\mathbb{Q}$ in $p^{-s}$.
\end{theorem}

In order to establish this, du Sautoy showed that the `analytic' theory of $p
$-adic analytic groups can be reduced to `$p$-adic analytic' model theory, as
developed by Denef and van den Dries. As well as opening up a fascinating new
field of study, this result led the way to some remarkable applications in the
theory of finite $p$-groups, discussed in the following section.

The study of these group-theoretic zeta functions is a very active area of
research at the present time; many results have been obtained but many more
challenging problems remain open. For more details and references up to 2002
see [dSS] and [SG], Chapters 15 and 16 (but there has been much progress since
then, for example in the work of Christopher Voll [V]).

Instead of counting subgroups of finite index, one could count equivalence
classes of finite-dimensional representations; the Dirichlet series encoding
these numbers give rise to `representation-growth zeta functions'. See [J-Z1],
[LL] and [V] for some recent and deep results about these.

\section{Finite $p$-groups\label{p-groups}}

\subsection{Coclass}

It was clear from the early days of group theory that the finite simple groups
are rather special: they are, essentially, the symmetry groups of highly
symmetrical structures (a finite set, or a vector space with a bilinear form).
Of course this wasn't actually proved until the 1980s (and the final steps
have only just been published), but the fact is that these objects form an
elegant list of identifiable objects, and they are `rigid' in two senses: (1)
they are \emph{isolated}: you can't move from one to the next by a small
deformation, and (2) the possibilities of building composite groups out of
them are very limited: they have small Schur multipliers and small outer
automorphism groups.

Neither of these (slightly vague) statements is true of nilpotent groups. It
was equally clear, at least from the 1930s with the work of Philip Hall and
others, that the finite $p$-groups constitute a vast and rather amorphous
collection. Thus the received wisdom for most of the last century considered
finite $p$-groups to be unclassifiable.

This pessimistic conclusion was based on the experience of trying to produce
coherent lists of $p$-groups, starting with the smallest and working up by
size; in practice this was only achieved for groups of nilpotency class 2 and
quite modest size, as the number of groups of order $p^{n}$ was found to grow
extremely fast with $n$. Higman and Sims showed in the 1960s that this number
is about $p^{\frac{2}{27}n^{3}}$, and that the number of groups of class 2 is
already about this big. (Contrast this with the number of simple groups of
order $n$, which is nearly always zero, sometimes one and very occasionally two!)

However, a different picture appears if instead of small nilpotency class one
looks at $p$-groups of \emph{large} class. Completing earlier work of
Blackburn, Leedham-Green and McKay found that the $p$\emph{-groups of maximal
class} do form a comprehensible pattern, and can indeed be neatly classified
by their order. What emerged from this classification is that, for a fixed
prime $p$, the best way to think of $p$-groups of maximal class is as the
finite quotients of one particular pro-$p$ group; for example, the $2$-groups
of maximal class are precisely the finite quotients of the `dihedral pro-$2$
group' $\mathbb{Z}_{2}\rtimes C_{2},$ together with certain natural
`twistings' of them (quaternion or semi-dihedral groups). This realisation led
Leedham-Green and Newman to formulate an audacious generalization, that became
known as the ``coclass conjectures''. These profoundly insightful conjectures
cast the problem of classifying $p$-groups into a completely new framework,
and totally transformed the subject between 1980 and 1994, when the
conjectures were finally established.

A finite $p$-group is said to have \emph{coclass} $r$ if it has order $p^{n}$
and nilpotency class $n-r$ (so maximal class means coclass $1$). A pro-$p$
group has coclass $r$ if it is the inverse limit of a system of finite
$p$-groups of coclass $r$ (with all maps surjective). The main conjecture of
Leedham-Green and Newman, Conjecture A, is purely finitary: it states that
every $p$-group of coclass $r$ has a normal subgroup of nilpotency class at
most $2$ and bounded index (the bound depending only on $p$ and $r$). In view
of the remarks above, this might seem like no progress as regards the
classification: what lies behind it, however, is a vision of the whole
universe of $p$-groups of fixed coclass. For given $p$ and $r,$ one arranges
the set of all (isomorphism types) of coclass $r$ $p$-groups into a graph
$\mathcal{G}(p,r)$, whose directed edges represent the quotient maps
$G\rightarrow G/Z$ where $Z$ is a central subgroup of order $p$ in $G$. Each
infinite chain in this graph then gives rise in a natural way to a\ pro-$p$
group of coclass $r$. Now the remarkable facts are these:

\begin{itemize}
\item There are only finitely many infinite pro-$p$ groups of coclass $r$ (for
given $p$ and $r$);

\item Each infinite pro-$p$ group of finite coclass is finitely generated and
virtually abelian, in other words, it is a finite extension of $\mathbb{Z}%
_{p}^{d}$ for some finite $d$.
\end{itemize}

\noindent Moreover, every finite $p$-group of coclass $r$ is either a quotient
of one of these virtually abelian pro-$p$ groups, or is obtained from such a
quotient by an explicit `twisting' process, or is one of finitely many
`sporadic' groups.

A key step in the proof, achieved by Leedham-Green, was to show that every
pro-$p$ group of finite coclass is a $p$\emph{-adic analytic group}, that is,
a pro-$p$ group of finite rank. Once this was known, it became possible to
apply powerful techniques for studying such groups, to show that if a $p$-adic
analytic group has finite coclass then it must be virtually abelian. The first
proof of this fact, due to Donkin, rests on the `analytic' aspect of these
groups and applies the classification of semisimple $p$-adic Lie algebras,
thus establishing a bridge between the theory of $p$-groups and the theory of
finite simple groups. Subsequently, a clever direct argument (also using Lie
algebras) was found by Shalev and Zelmanov. An alternative, purely finitary,
proof for Conjecture A was later obtained by Shalev [Sh]; although this avoids
the use of pro-$p$ groups altogether, it was clearly inspired by the $p$-adic
methods used before.

Explanatory accounts of all or parts of this story are to be found in [LGM1],
[LGM2], [DDMS], Chapter 10. For full references to the many original papers,
see the bibliographies to [LGM1] and [LGM2].

\subsection{Conjecture P}

The main results of coclass theory show that the graph $\mathcal{G}(p,r)$ has
finitely many components; moreover, if we remove a finite number of `sporadic'
groups what remains is the disjoint union of finitely many trees. Each of
these trees contains just one maximal infinite chain, the `trunk', to which
are attached infinitely many finite `twigs'. On the basis of extensive
computer investigations, M. Newman and E. O'Brien were led to make some very
precise conjectures about the shape of these trees. In particular, their
\emph{Conjecture P} asserts that when $p=2$, each tree is eventually periodic,
with period dividing $2^{r-1}$.

The conjecture obviously implies that the twigs of such a tree are of bounded
length, and this is no longer true when the prime $p$ is odd. However, du
Sautoy was able to establish a general periodicity result which includes (the
qualitative part of) Conjecture P as a special case. For each tree
$\mathcal{T}$ as above and each natural number $m$, let $\mathcal{T}[m]$
denote the `pruned tree' obtained from $\mathcal{T}$ by removing all vertices
whose distance from the trunk exceeds $m$.

\begin{theorem}
\emph{[dS2]} Each of the pruned trees $\mathcal{T}[m]$ is eventually periodic.
\end{theorem}

\noindent It is known that when $p=2$ the twigs have bounded lengths, so in
this case we have $\mathcal{T}[m]=\mathcal{T}$ for some value of $m$.

The proof is a remarkable application of du Sautoy's rationality theorem for
zeta functions (see \S \ref{sg} above). First of all, he deduces from the
results of coclass theory that there exists a certain $p$-adic analytic
pro-$p$ group $H=H(p,r)$ which maps onto every finite $p$-group of coclass
$r$. The holomorph $P=H\rtimes\mathrm{Aut}(H)$ is again a $p$-adic analytic
group, and du Sautoy associates a certain generalized zeta function to the
pair $(H,P)$; the coefficients of (the Dirichlet series defining) this
function encode precisely the `shape' of the pruned tree $\mathcal{T}[m]$. He
proves that this generalized zeta function is again rational, and the stated
periodicity then emerges as a formal consequence. For details of this
argument, see [dSS].

Zeta functions are also used in [dS2] to obtain results about the enumeration
of $p$-groups and of finite nilpotent groups of fixed \emph{nilpotency class}.
These are also discussed in [dSS].

\section{Finitely generated groups\label{fg}}

\subsection{Linearity}

Nearly a century ago, Hasse argued in favour of treating the $p$-adic
completions of $\mathbb{Q}$ on the same footing as the reals. This idea had a
huge influence on the development of number theory; as mentioned above, it
also led to the idea of studying (non-commutative) groups via their profinite
completions. In general, a group deosn't even have a `real completion' (unless
it is nilpotent, say), but every group has its pro-$p$ completions and its
profinite completion. Thus every group can be mapped, functorially, into
various interesting compact topological groups.

This simple idea led Lubotzky to the solution of a long-standing problem: how
to characterize, by purely internal criteria, those groups that have a
faithful finite-dimensional linear representation over some field, `linear
groups' for short.

\begin{theorem}
\emph{[Lu2] `Lubotzky linearity criterion'} Let $G$ be a finitely generated
group. Then $G$ is linear over some field of characteristic zero if and only
if, \ for some prime $p$ and some natural number $r$, $G$ has a chain of
normal subgroups
\[
G=G_{0}\geq G_{1}\geq G_{2}\geq\ldots
\]
such that \emph{(i) }$G/G_{1}$ is finite, \emph{(ii)} $G_{1}/G_{n}$ is a
finite $p$-group of rank at most $r$ for every $n\geq1$, and \emph{(iii)}
$\bigcap_{n}G_{n}=1$.
\end{theorem}

Suppose that $G$ satisfies the given condition, and consider the inverse
limit
\[
P=\underset{\longleftarrow}{\lim}G_{1}/G_{n}.
\]
Hypothesis (ii) implies that this group $P$ is a pro-$p$ group of finite rank,
and so a $p$-adic analytic group (see \S \ref{pag} above). Then Lie theory and
Ado's Theorem show that $P$ is linear over the $p$-adic number field
$\mathbb{Q}_{p}$. Hypothesis (iii) implies that $G_{1}$ embeds into $P$, so
$G_{1}$ is linear, and it follows by hypothesis (i) that $G$ itself is linear
(form the induced representation).

Note that the argument so far does not require $G$ to be finitely generated;
the converse, however does. To see why it is true, suppose now that $G$ is a
finitely generated subgroup of $\mathrm{GL}_{d}(F)$ where $F$ is a field of
characteristic zero. Then in fact $G\leq\mathrm{GL}_{d}(R)$ where $R$ is some
finitely generated subring of $F$. Commutative algebra shows that for almost
all primes $p$, such a ring $R$ can be embedded in a matrix ring over
$\mathbb{Z}_{p}$; for each such prime it follows that $G$ can be embedded in
some $\mathrm{GL}_{d^{\prime}}(\mathbb{Z}_{p})$ (where $d^{\prime}=md$ may
depend on $p$). Choosing a suitable prime $p$ and identifying $G$ with its
image in $\mathrm{GL}_{d^{\prime}}(\mathbb{Z}_{p})$, we take
\[
G_{n}=\left\{  g\in G\mid g\equiv\mathbf{1}\,(\operatorname{mod}%
p^{n})\right\}  .
\]
It is easy to see that the sequence $(G_{n})$ then satisfies conditions (i)
and (iii); and condition (ii) is satisfied because the `first congruence
subgroup'
\[
\mathrm{GL}_{d^{\prime}}^{1}(\mathbb{Z}_{p})=\ker\left(  \mathrm{GL}%
_{d^{\prime}}(\mathbb{Z}_{p})\rightarrow\mathrm{GL}_{d^{\prime}}%
(\mathbb{Z}_{p}/p\mathbb{Z}_{p})\right)
\]
is a pro-$p$ group of finite rank ([DDMS], Chapter 5).

For a more detailed account, and several variations on the same theme, see
[DDMS] Interlude B.

So far, no-one has succeeded in establishing, or even formulating, an
analogous characterization of the finitely generated linear groups over fields
of positive characteristic, and this remains a challenging open problem.
Lubotzky's criterion can paraphrased as: ``some pro-$p$ completion of some
normal subgroup of finite index in $G$ is $p$-adic analytic''; a natural
starting point for the characteristic-$p$ analogue would be to gain a better
understanding of the pro-$p$ groups that are `analytic' over a local ring of
characteristic $p$; the beginnings of such a theory are outlined in the final
chapter of [DDMS].

The `Lubotzky criterion' arises from considering congruence subgroups modulo
powers of a fixed prime -- looking `downwards', we may say. Another way of
looking at a finitely generated linear group is `sideways': for example, we
can embed $\mathrm{GL}_{d}(\mathbb{Z})$ into the Cartesian product $\prod
_{p}\mathrm{GL}_{d}(\mathbb{Z}/p\mathbb{Z})$, where $p$ ranges over any
infinite set of primes. A. I. Mal'cev generalized this observation to show
that every finitely generated linear group of degree $d$ is residually `linear
of degree $d$ over a finite field'. The precise converse is not true, but J.
S. Wilson showed that a slightly weaker statement does hold: \emph{if a
finitely generated group }$G$ \emph{is residually (linear of degree}
$d$\emph{) then} $G$ \emph{is a subdirect product of finitely many linear
groups}. For the proof, and some refinements, see [SG], Window 8. Here is one
such refinement, which serves as a reduction step for many of the results
stated in \S \ref{sg}, above:

\begin{theorem}
\label{reslin}Let $G$ be a finitely generated group and $(N_{i})$ a family of
normal subgroups of $G$ with $\bigcap_{i}N_{i}=1$. Suppose that $G/N_{i}%
\leq\mathrm{GL}_{d}(F_{i})$, where each $F_{i}$ is either a field of
characteristic zero or a finite field, and suppose further that for each prime
$p$ the number of $i$ with $\mathrm{char}F_{i}=p$ is finite. Then $G$ is
linear over a field of characteristic zero.
\end{theorem}

\subsection{Finite quotients\label{fq}}

What does it mean for a family of finite groups $\mathcal{X}$ to be precisely
the set $\mathcal{F}(\Gamma)$ of (isomorphism types of) all finite quotients
of some finitely generated group $\Gamma$? Equivalently, what does it mean for
a profinite group $G$ to be the profinite completion of a finitely generated
(abstract) group? As mentioned in \S \ref{local}, this holds if and only if
$G$ contains a dense finitely generated subgroup $\Gamma$ that has the
congruence subgroup property; so the question may be seen as finding necessary
and/or sufficient conditions on a profinite group $G,$ expressed in terms of
the family $\mathcal{F}(G)$, for the existence of such a subgroup (when $G$ is
a \emph{profinite} group, $\mathcal{F}(G)$ denotes the set of
\emph{continuous} finite quotient groups of $G$: see the next subsection).

Two obvious necessary conditions for such a family $\mathcal{X}$ are (1) that
$\mathcal{X}$ is quotient-closed, and (2) that all the groups in $\mathcal{X}$
can be generated by some bounded number of elements; but it seems very
difficult to find further, less obvious ones. Suppose for example that
$\mathcal{X}$ contains a subgroup $X_{i}$ of $\mathrm{GL}_{d}(F_{i})$ for
$i=1,2,\ldots$ where $F_{i}$ is a finite field of characteristic $p_{i}$ and
$p_{1},p_{2},\ldots$ \ is an infinite sequence of distinct primes. Then
$\Gamma$ has a quotient $\overline{\Gamma}$ which satisfies the hypotheses of
Theorem \ref{reslin}, so $\overline{\Gamma}$ is a finitely generated
characteristic-zero linear group; if we assume also that the groups $X_{i}$
are simple and of unbounded orders (or some suitable weaker condition), we
find that $\overline{\Gamma}$ is not virtually soluble. Under these
conditions, $\overline{\Gamma}$ is guaranteed to possess a host of special
finite quotients: applying a deep `strong approximation' theorem due to Nori
and Weisfeiler, Lubotzky established the following important result:

\begin{theorem}
\label{LA}\emph{`Lubotzky alternative'} Let $\Gamma$ be a finitely generated
linear group over a field of characteristic zero. Then one of the following holds:

\emph{(a)} $\Gamma$ is virtually soluble;

\emph{(b)} there exist a connected, simply connected simple algebraic group
$\mathfrak{G}$ over $\mathbb{Q}$, a finite set of primes $S$ such that
$\mathfrak{G}(\mathbb{Z}_{S})$ is infinite, and a subgroup $\Gamma_{1}$ of
finite index in $\Gamma$ such that the profinite group $\mathfrak{G}%
(\widehat{\mathbb{Z}_{S}})$ is an image of $\widehat{\Gamma_{1}}$.
\end{theorem}

\noindent(Here $\mathbb{Z}_{S}=\mathbb{Z}[\frac{1}{p};\,p\in S]$, and
$\mathfrak{G}(\widehat{\mathbb{Z}_{S}})$ is isomorphic to the product
$\prod_{p\notin S}\mathfrak{G}(\mathbb{Z}_{p})$.) For the proof, see [SG],
Window 9. Applying this to the group $\overline{\Gamma}$, we may deduce that
the set $\mathcal{X}$ must contain many other groups in addition to the
$X_{i}$: for each prime $p\notin S$ and each $n$, a group $Q_{p}$ containing
$\mathfrak{G}(\mathbb{Z}/p^{n}\mathbb{Z})$ as a subgroup, the indices
$\left\vert Q_{p}:\mathfrak{G}(\mathbb{Z}/p^{n}\mathbb{Z})\right\vert $ being
bounded above by a constant.

Thus if $\mathcal{P}$ is an infinite set of primes, a set of groups like%
\[
\left\{  \prod_{p\in T}\mathrm{PSL}_{d}(\mathbb{F}_{p})\mid T\text{ a finite
subset of }\mathcal{P}\right\}
\]
cannot be the whole of $\mathcal{F}(\Gamma)$ for a finitely generated group
$\Gamma$, while of course it is equal to $\mathcal{F}(G)$ where $G=\prod
_{p\in\mathcal{P}}\mathrm{PSL}_{d}(\mathbb{F}_{p})$. Thus the $2$-generator
profinite group $G$ cannot be the profinite completion of a finitely generated group.

The problem with this group $G$ is that the finite simple factor groups have
bounded ranks. In an amazing feat of ingenuity, Kassabov and Nikolov have
recently shown that this is essentially the \emph{only} obstacle, when it
comes to products of finite simple groups. For a group $S$ they write $l(S)$
to denote the largest integer $k$ such that $S$ contains a copy of the
alternating group $\mathrm{Alt}(k)$, and they prove

\begin{theorem}
\label{KN}\emph{[KN] }Let $(S_{n})$ be a sequence of finite simple groups such
that $l(S_{n})\rightarrow\infty$ as $n\rightarrow\infty$, and let%
\[
G=\prod_{n=1}^{\infty}S_{n}\text{.}%
\]
If $G$ is finitely generated (as a profinite group), then $G$ is the profinite
completion of a finitely generated group.
\end{theorem}

Thus we now have some understanding of what it means for a family of
\emph{semisimple} groups (products of finite simple groups) to be equal to
$\mathcal{F}(\Gamma)$ for some finitely generated group $\Gamma$. However, it
seems difficult even to formulate a conjecture about the nature of sets like
$\mathcal{F}(\Gamma)$ in general.

Meanwhile, we could consider weakening the question a little, and asking: what
does it mean for a collection of finite simple groups to be precisely the
collection of \emph{composition factors }of groups in $\mathcal{F}(\Gamma)$
for some finitely generated group $\Gamma$? These are called the \emph{upper
composition factors} of $\Gamma$. An almost complete answer is provided by

\begin{theorem}
\emph{[S4] }Let $\mathcal{S}$ be any collection of (isomorphism types) of
non-abelian finite simple groups. Then there exists a $63$-generator group
$\Gamma$ whose set of upper composition factors is precisely $\mathcal{S}.$
\end{theorem}

To construct such a group $\Gamma$ we start with a suitable profinite group
$G$, and then find $\Gamma$ as a dense subgroup in $G$. To begin with, we
enumerate $\mathcal{S}$ as $\{X_{1},X_{2},\ldots,X_{n},\ldots\}$ (if
$\mathcal{S}$ is finite, the result is trivial, given that every finite simple
group can be generated by $2$ elements: this follows from CFSG, and implies
that $\Gamma=\prod_{X\in\mathcal{S}}X$ is a $2$-generator group). For each $n$
we pick a faithful primitive permutation representation for $X_{n}$, and so
identify $X_{n}$ with a subgroup of $\mathrm{Sym}(l_{n})$ for some $l_{n}$.
Take $W_{1}=X_{1},$ for $n>1$ let $W_{n}$ be the permutational wreath product
\[
W_{n}=X_{n}\wr W_{n-1},
\]
and define
\[
G=\underset{\longleftarrow}{\lim}W_{n}.
\]
Thus $G$ is a profinite group, whose set of upper composition factors is
precisely $\mathcal{S}$.

This is all very easy. The challenge now is to find a suitable dense subgroup
in $G$. The key lies in realizing $G$ as a group of automorphisms of a
suitable object.

Given the sequence of positive integers $(l_{n})$, consider the
\emph{spherically homogeneous rooted tree} $\mathcal{T}$ of type $(l_{n})$:
this is a connected graph without circuits, having a distinguished vertex
$v_{0}$ (the \emph{root}), and for each $n\geq1$ having $l_{1}\ldots l_{n}$
vertices at distance $n$ from the root, each of valency $1+l_{n+1}$ (so at
each vertex of `level' $n\geq1$ there is one edge pointing `upwards' towards
the root and $l_{n+1}$ edges pointing `downwards' to the next level). It is
easy to see that the automorphism group of this structure is the inverse limit
of the finite permutational wreath products
\[
V_{n}=\mathrm{Sym}(l_{n})\wr\ldots\wr\mathrm{Sym}(l_{2})\wr\mathrm{Sym}%
(l_{1}).
\]
Thus $V_{n}$ contains $W_{n}$ as a permutation group for each $n$, and we may
identify our profinite group $G$ as a closed subgroup of $\mathrm{Aut}%
(\mathcal{T})$; a base for the neighbourhoods of the identity in $G$ is given
by the `level-stabilizers' $\mathrm{st}_{G}(n)=\ker\left(  G\rightarrow
W_{n}\right)  $.

One of the main results of [S4] states that there exists a $61$-generator
perfect group $P$ that maps onto every non-abelian finite simple group. Using
this, we define $63$ specific tree automorphisms of $\mathcal{T}$, all lying
in the group $G$, and take $\Gamma$ to be the group generated by these $63$
automorphisms. These generators are so chosen that (a) for each $n $, the
group $\Gamma$ acts as the whole group $W_{n}$ on the $n$th level of
$\mathcal{T},$ and (b) each nontrivial normal subgroup of $\Gamma$ contains
$\mathrm{st}_{\Gamma}(n)=\Gamma\cap\mathrm{st}_{G}(n)$ for some $n$ (actually,
a quite general argument shows that each nontrivial normal subgroup of
$\Gamma$ contains the derived group $\mathrm{st}_{\Gamma}(n)^{\prime}$ of
$\mathrm{st}_{\Gamma}(n)$ for some $n$; the role of the perfect group $P$ is
to ensure that in our case we have $\mathrm{st}_{\Gamma}(n)^{\prime
}=\mathrm{st}_{\Gamma}(n)$ for each $n$). Property (a) means that $\Gamma$ is
dense in $G$, while property (b) implies that $\Gamma$ has the CSP in $G$. It
follows that $\mathcal{F}(\Gamma)=\mathcal{F}(G)$, and hence that the set of
upper composition factors of $\Gamma$ is precisely $\mathcal{S}$.

The same construction, using sets $\mathcal{S}$ of the form $\{\mathrm{PSL}%
_{2}(\mathbb{F}_{p})\mid p\in\mathcal{P}\}$ for suitably chosen sets of primes
$\mathcal{P}$, was used in [S4] to construct finitely generated groups with
arbitrarily specified types of subgroup growth (within a certain range). For
details, and more discussion of trees like $\mathcal{T}$, see Chapter 13 of [SG].

Certain groups of rooted tree automorphisms called \emph{branch groups} have
been studied in depth by Grigorchuk and others. These include the groups
described above, but are more usually pro-$p$ groups (or dense finitely
generated subgroups thereof); the celebrated construction by Grigorchuk of a
finitely generated group having `intermediate word growth' was (a dense
subgroup of) a pro-$2$ branch group. See [G] and [BG].

\subsection{Forgetting the topology}

To be given a profinite group $G$ is more or less equivalent to being given
the family of all finite continuous quotient groups of $G$, that is, the
groups $G/N$ where $N$ ranges over all the open normal subgroups of $G.$
Indeed, $G$ is (naturally isomorphic to) the inverse limit of this family,
relative to the natural quotient maps $G/N\rightarrow G/M$, \ $(M\geq N)$. If
we forget the topology and think of $G$ just as an abstract group, we would
expect to lose a lot of information: out of all the normal subgroups of finite
index in $G$, how could we possibly pick out those that were open? Consider
the following simple example. Fix a prime $p$, for each $i$ let $C_{i}$ be
cyclic of order $p$, put $G_{n}=C_{1}\times\cdots\times C_{n}$ and let
\[
G=\underset{\longleftarrow}{\lim}G_{n}%
\]
where $G_{m}\rightarrow G_{n}$ for $m\geq n$ are the obvious projection maps.
The open subgroups of $G$ are those that contain $\ker(G\rightarrow G_{n})$
for some $n$, so there are only countably many of them. On the other hand, as
an abstract group $G$ is abelian, of exponent $p$ and uncountable (of
cardinality $\mathfrak{c}=2^{\aleph_{0}}$); it is therefore a $\mathfrak{c}%
$-dimensional vector space over $\mathbb{F}_{p}$ and so contains
$2^{\mathfrak{c}}$ subspaces of finite codimension. Thus $G$ has
$2^{\mathfrak{c}}$ (normal) subgroups of finite index, of which only countably
many are open. It is obvious, from the very homogeneous nature of (the
abstract group) $G$, that there is no way of recovering the original topology.
(A similar construction can be made using any nontrivial finite group in place
of the group of order $p$: see [RZ2], Ex. 4.2.13.)

However: if we restrict attention to (topologically) \emph{finitely generated
}profinite groups, the opposite is true:

\begin{theorem}
\label{fip}\emph{[NS2] }In a finitely generated profinite group, every
subgroup of finite index is open.
\end{theorem}

\noindent This is a remarkable fact: if we form the inverse limit $G$ of any
(surjective) inverse system $\mathcal{S}$ of finite groups, all of which can
be generated by some fixed number of elements, then the only finite groups
onto which $G$ can be mapped homomorphically are the quotients of groups in
$\mathcal{S}$; moreover, since the subgroups of finite index form a base for
the neighbourhoods of the identity, the topology of $G$ is completely
determined by its structure as an abstract group.

This theorem is a case where a problem on profinite groups served as the
motivation for some new developments in finite group theory, and it
illustrates very clearly the principle that a qualitative property of
profinite groups corresponds to a \emph{uniform} quantitative property of
finite groups. The basic idea is as follows. Let $w=w(x_{1},\ldots,x_{k})$ be
a group word, and $G$ a profinite group. Since the mappings
\begin{align*}
(g_{1},\ldots,g_{k})  &  \mapsto w(g_{1},\ldots,g_{k}),\\
(g_{1},\ldots,g_{k})  &  \mapsto w(g_{1},\ldots,g_{k})^{-1}%
\end{align*}
from $G^{(k)}$ to $G$ are continuous, their images in $G$ are compact. It
follows that for each $n$, the set $S(n)$ of all products of $n$ elements of
the form $w(g_{1},\ldots,g_{k})^{\pm1}$ is compact, hence closed in $G$. Now
consider the \emph{verbal subgroup} $w(G)$, generated (\emph{algebraically},
not topologically) by all values of $w$ in $G$:
\begin{equation}
w(G)=\bigcup_{n=1}^{\infty}S(n). \label{union}%
\end{equation}
If it happens that for some finite $n$ we have $w(G)=S(n)$, then $w(G)$ is
closed; conversely, if $w(G)$ is closed then a simple argument using the Baire
category theorem and (\ref{union}) shows that $w(G)=S(n)$ for some $n$. This
means that every product of $w$-values in $G$ is equal to a product of $n$
$w$-values (where by `$w$-value' I mean an element of the form $w(g_{1}%
,\ldots,g_{k})^{\pm1}$); let me abbreviate this to `$w$ has width $n$ in $G$'.

On the other hand, it is easy to see that $w$ has width $n$ in $G$ if and only
if $w$ has width $n$ in $G/N$ for every open normal subgroup $N$ of $G$.
Indeed, if the latter holds then
\[
\frac{w(G)N}{N}=w(G/N)=\frac{S(n)N}{N}%
\]
for each $N$, so
\begin{equation}
S(n)\subseteq w(G)\subseteq\bigcap_{N}w(G)N=\bigcap_{N}S(n)N=S(n),
\label{intersection}%
\end{equation}
the last equality holding because $S(n)$ is a closed subset of $G$. The
converse is obvious. Thus we have established the link between a qualitative
property of $G$ and a uniform property of $\mathcal{F}(G)$ (the set of
\emph{continuous} finite images of $G$):

\begin{proposition}
\label{baire}Let $G$ be a profinite group and $w$ a group word. Then the
(algebraic) verbal subgroup $w(G)$ is closed in $G$ if and only if there
exists a natural number $n$ such that $w$ has width $n$ in every continuous
finite image of $G$.
\end{proposition}

This result, due to Brian Hartley, is nice, but how does it help with our
original problem? Suppose we know in addition that the index $\left|
Q:w(Q)\right|  $ is uniformly bounded for all $Q\in\mathcal{F}(G)$. Then the
big intersections in the middle of (\ref{intersection}) contain only finitely
many distinct terms, each of which is an open subgroup of $G$; and we may
infer that in this case, $w(G)$ is not only closed but \emph{open}.

Now let $G$ be a $d$-generator profinite group and $H$ a subgroup of finite
index. Then $H$ contains a subgroup $K$ which is normal and of finite index in
$G$. Let $F=F_{d}$ be the free group on free generators $x_{1},\ldots,x_{d}$
and let $D$ be the intersection of the kernels of all homomorphisms
$F\rightarrow G/K.$ Then $D$ has finite index in $F$ and is therefore finitely
generated, by $w_{1}(x_{1},\ldots,x_{d}),\ldots,w_{m}(x_{1},\ldots,x_{d})$
say. Put
\[
w(\mathbf{y}_{1},\ldots,\mathbf{y}_{m})=w_{1}(\mathbf{y}_{1})\ldots
w_{m}(\mathbf{y}_{m})
\]
where $\mathbf{y}_{1},\ldots,\mathbf{y}_{m}$ are disjoint $d$-tuples of
variables. It is easy to see that (i) $w(F)=D$ and (ii) $w(G)\leq K$. The
latter shows that $H$ will be open in $G$ if $w(G)$ is open. Property (i)
implies that
\[
\left|  Q:w(Q)\right|  \leq\left|  F:w(F)\right|  <\infty
\]
for every $Q\in\mathcal{F}(G)$ (while we know nothing, a priori, about the
finite group $G/K$, we do know that each of the finite groups in
$\mathcal{F}(G)$ is $d$-generator, hence an image of $F$). To conclude that
$w(G)$, and therefore also $H$, is open in $G$, we are thus reduced to
establishing the following `uniformity theorem' about finite groups (I call
$w$ `$d$-locally finite' if $\left|  F_{d}:w(F_{d})\right|  $ is finite):

\begin{theorem}
\label{unif}\emph{[NS2]} Let $d$ be a natural number and let $w$ be a
$d$-locally finite group word. Then there exists $f=f(w,d)$ such that $w$ has
width $f$ in every $d$-generator finite group.
\end{theorem}

The proof of this result is long and difficult, and depends on CFSG. I will
say no more about it here; for a brief outline see the announcement [NS1].

In the same paper we establish an analogous theorem for the commutator words
$w=[x_{1},\ldots,x_{k}]$; in view of Proposition \ref{baire} this implies that
the derived group, and the higher terms of the lower central series, are
closed in every finitely generated profinite group. We also made a not
entirely successful attempt to do the same for the words $w=x^{q}$
($q\in\mathbb{N}$), so the following is still open:\medskip

\noindent\textbf{Problem }Let $q$ be a natural number. Is it true that the
subgroup $G^{q}=\left\langle g^{q}\mid g\in G\right\rangle $ (generated
\emph{algebraically} by all $q$th powers in $G$) is open in $G$, for every
(topologically) finitely generated profinite group $G$?\medskip

\noindent Note that in this situation, $G^{q}$ is open if and only if it is
closed, because there is a finite upper bound for the order of every finite
$d$-generator group of exponent dividing $q$: this is the positive solution of
the restricted Burnside Problem, due to Zelmanov. Thus the problem is
equivalent to asking whether, for each $d$, the \textquotedblleft Burnside
word\textquotedblright\ $x^{q}$ has bounded width in all $d$-generator finite groups.

Whatever the answer turns out to be, results of this type certainly don't hold
for arbitrary words: Romankov [R] has given a simple construction for a
three-generator soluble pro-$p$ group $G$ in which the second derived group
$G^{\prime\prime}$ is not closed; and $G^{\prime\prime}=w(G)$ where
$w=[[x_{1},x_{2}],[x_{3},x_{4}]]$. It would be very interesting to find a
characterization of those group words $w$ which have the uniformity property
of Theorem \ref{unif}. This is equivalent to asking: \emph{for which words
}$w$\emph{ is it the case that }$w(G)$\emph{ is closed in }$G$\emph{ for every
finitely generated profinite group }$G$ \emph{?}; if we restrict to pro-$p$
groups, the remarkably simple answer has recently been discovered by Andrei
Jaikin [J-Z2]: $w(G)$ is closed in $G$ for every finitely generated pro-$p$
group $G$ if and only if $w\notin F^{\prime\prime}(F^{\prime})^{p}$, where $F$
is the free group on the variables occurring in $w$.

Let us turn briefly to the non-finitely generated case. For a profinite group
$G,$ let $G_{0}$ denote the underlying abstract group. Theorem \ref{fip}
implies that if $G$ is finitely generated then $\mathcal{F}(G)=\mathcal{F}%
(G_{0})$ (recall that these denote the sets of \emph{isomorphism types} of
finite quotients, by open normal subgroups or by all normal subgroups of
finite index, respectively). We have also seen examples of (infinitely
generated) profinite groups $G$ that have many non-open normal subgroups of
finite index; but in these examples, too, we have $\mathcal{F}(G)=\mathcal{F}%
(G_{0})$ -- the same finite groups appear, though with different
multiplicities as quotients of $G$. To construct a group $G$ such that
$\mathcal{F}(G)\neq\mathcal{F}(G_{0})$ takes a little more effort; the
following example was suggested by Lubotzky and Holt. For a finite group
$S=S^{2}$, let $f(S)$ denote the least integer $n$ such that every element of
$S$ is equal to a product of $n$ squares (here $S^{2}$ denotes the subgroup
generated by all squares). Now for each $n$ let $S_{n}$ be a finite group with
$S_{n}=S_{n}^{2}$ and $f(S_{n})>n$, and take $G=\prod_{n=1}^{\infty}S_{n}$.
Proposition \ref{baire} shows that the subgroup $G^{2}$ is not closed in $G$;
in particular it can't be equal to $G$, so $G/G^{2}$ has the cyclic group
$C_{2}$ of order $2$ as a quotient. On the other hand, $C_{2}\notin
\mathcal{F}(G)$ since every continuous finite quotient of $G$ is a quotient of
$S_{1}\times\cdots\times S_{k}$ for some $k$. Thus $C_{2}\in\mathcal{F}%
(G_{0})\setminus\mathcal{F}(G)$.

Suitable groups $S_{n}$ may be constructed as follows (for details, see [H]).
Let $H=\mathrm{SL}_{2}(\mathbb{F}_{4})$ and let $M$ be its natural
$2$-dimensional $\mathbb{F}_{4}$-module, considered as a $4$-dimensional
$\mathbb{F}_{2}H$-module. There is an $H$-epimorphism $\phi$ from
$M\otimes_{\mathbb{F}_{2}}M$ onto the trivial module $\mathbb{F}_{2}$. Now let
$M_{1},\ldots,M_{k}$ be copies of $M$ and form a special $2$-group $P$ with
$P/[P,P]=M_{1}\times\cdots\times M_{k}$ and $\mathrm{Z}(P)=[P,P]=\prod
_{i<j}[M_{i},M_{j}]$, where $[M_{i},M_{j}]\cong\mathbb{F}_{2}$ and the
commutator mapping $M_{i}\times M_{j}\rightarrow\lbrack M_{i},M_{j}]$ for
$i<j$ is induced by $\phi$. Then $H$ acts by automorphisms on $P$, fixing
$\mathrm{Z}(P)$ elementwise, and we set $S_{n}=P\rtimes H$. It is easy to see
that $S_{n}=[S_{n},S_{n}]=S_{n}^{2}$. Since $(zx)^{2}=x^{2}$ for every
$z\in\mathrm{Z}(P)$ and $x\in S_{n}$, the number of squares in $S_{n}$ is no
more than $\left|  S_{n}/\mathrm{Z}(P)\right|  =4^{k}\cdot60$; on the other
hand $\left|  S_{n}\right|  =2^{k(k-1)/2}\cdot4^{k}\cdot60$. This implies that
$f(S_{n})\geq(k+3)/6>n$ if we choose $k>6n$.

Let me conclude with a little exercise for the reader: if $G$ is any profinite
group, then every group in $\mathcal{F}(G_{0})$ is isomorphic to a
\emph{section} of some group in $\mathcal{F}(G)$ (hint: apply Theorem
\ref{fip} to a suitable finitely generated subgroup of $G$).

\section{Probability\label{prob}}

Every compact topological group has an invariant measure, the \emph{Haar
measure}, unique up to a multiplicative constant. Though quite tricky to
construct in general, it is very easy to evaluate in the special case of a
profinite group $G$. Let us write $\mu(X)$ for the measure of a subset $X$ of
$G$, and normalize $\mu$ so that $\mu(G)=1$. If $H$ is an open subgroup of $G$
then each coset $Hx$ of $H$ has the same measure, so
\[
\mu(Hx)=\left|  G:H\right|  ^{-1}\mu(G)=\left|  G:H\right|  ^{-1}.
\]
Similarly, $\mu(xH)=\left|  G:H\right|  ^{-1}$. As the cosets of open
subgroups form a base for the open sets in $G,$ this determines the measure of
every open set, and hence also of every closed set. Assuming that $G$ is
countably based (i.e. has only countably many open normal subgroups) it is
easy to deduce that for any closed subset $X$ of $G$ we have
\begin{equation}
\mu(X)=\inf\frac{\left|  \pi(X)\right|  }{\left|  \pi(G)\right|  } \label{mu}%
\end{equation}
where $\pi$ ranges over all the quotient maps $G\rightarrow G/N$, $N$ an open
normal subgroup.

Now a measure space of measure $1$ is a \emph{probability space}: we interpret
$\mu(X)$ as the probability that a random element of $G$ belongs to the subset
$X$ (note that when $\pi(G)$ is finite, $\left|  \pi(X)\right|  /\left|
\pi(G)\right|  $ is just the proportion of elements of $\pi(G)$ that lie in
$\pi(X)$). So we can ask questions about the probability of interesting
group-theoretic events; for example, what is the probability that a random
$k$-tuple of elements generates $G$ (topologically)? To make this precise we
need to consider the measure on $G^{(k)}=G\times\cdots\times G$, still denoted
$\mu$, and define
\[
P(G,k)=\mu(X_{k})
\]
where
\[
X_{k}=\left\{  (x_{1},\ldots,x_{k})\in G^{(k)}\mid\left\langle x_{1}%
,\ldots,x_{k}\right\rangle =G\right\}  .
\]
(Here $\left\langle S\right\rangle $ denotes the closed subgroup of $G$
generated by the subset $S$.) The formula (\ref{mu}) becomes
\begin{equation}
P(G,k)=\inf P(G/N,k) \label{P}%
\end{equation}
where $N$ ranges over all open normal subgroups of $G$. Obviously, $P(G,k)=0$
unless $G$ can be generated by $k$ elements. But the converse is not always
true. Consider for example the procyclic group $G=\widehat{\mathbb{Z}}$, the
profinite completion of the infinite cyclic group $\mathbb{Z}$. Certainly $G$
can be generated by one element. On the other hand, it is easy to see that if
$n=p_{1}^{f_{1}}\ldots p_{r}^{f_{r}}$ then
\[
P(\mathbb{Z}/n\mathbb{Z},k)=\prod_{i=1}^{r}P(\mathbb{Z}/p_{i}^{f_{i}%
}\mathbb{Z},k)=\prod_{i=1}^{r}\left(  1-\frac{1}{p_{i}^{k}}\right)
\]
(since a subset $Y$ generates $\mathbb{Z}/p^{f}\mathbb{Z}$ unless $Y\subseteq
p\mathbb{Z}/p^{f}\mathbb{Z}$). Thus (\ref{P}) gives
\begin{align*}
P(\widehat{\mathbb{Z}},k)  &  =\prod_{p}\left(  1-\frac{1}{p^{k}}\right) \\
&  =\zeta(k)^{-1}\\
&  =\left\{
\begin{array}
[c]{ccc}%
0 &  & (k=1)\\
&  & \\
\frac{6}{\pi^{2}} &  & (k=2)
\end{array}
\right.  .
\end{align*}
The procyclic group $\widehat{\mathbb{Z}}$ is `only just' a one-generator
group: almost all elements do not generate it. On the other hand, a positive
proportion -- about $3/5$ -- of pairs do generate $\widehat{\mathbb{Z}}$.

Avinoam Mann calls a profinite group \emph{positively finitely generated}, or
PFG, if $P(G,k)>0$ for some natural number $k$. To get some feeling for this
property, note that $(x_{1},\ldots,x_{k})$ belongs to the set $X_{k}$ defined
above if and only if no maximal (open, proper) subgroup of $G$ contains all of
$x_{1},\ldots,x_{k}.$ That is,
\[
G^{(k)}\setminus X_{k}=\bigcup_{M\in\mathcal{M}}M^{(k)}
\]
where $\mathcal{M}$ denotes the set of all maximal subgroups of $G$. It
follows that
\begin{align*}
1-P(G,k)  &  =\mu\left(  \bigcup_{M\in\mathcal{M}}M^{(k)}\right) \\
&  \leq\sum_{M\in\mathcal{M}}\mu(M^{(k)})\\
&  =\sum_{M\in\mathcal{M}}\left|  G:M\right|  ^{-k}=\sum_{n\geq2}%
m_{n}(G)n^{-k}%
\end{align*}
where $m_{n}(G)$ is the \emph{number of maximal subgroups of index} $n$
\emph{in} $G$. Thus $P(G,k)$ is positive if the final sum is less than $1$.
Suppose for example that the numbers $m_{n}(G)$ grow at most like a power of
$n$ -- in this case $G$ is said to have \emph{polynomial maximal subgroup
growth}, or PMSG. Then for a certain $\alpha$ we have
\[
1-P(G,k)\leq\sum_{n\geq2}n^{\alpha-k}=\zeta(k-\alpha)-1
\]
which is less than $1$ if $k-\alpha\geq2$.

It follows that \emph{every profinite group with} PMSG \emph{is} PFG. Since
PMSG is a weaker condition than polynomial subgroup growth, we have the
corollary that \emph{every profinite group with} PSG \emph{is finitely
generated}. This fact can also be seen from the characterization of profinite
PSG groups, discussed in \S \ref{sg}, above; but it is remarkable that it
emerges from such a simple probabilistic argument. This simple argument is not
reversible, of course; a much more difficult argument, using detailed
information about the maximal subgroups of finite simple groups, enabled Mann
and Shalev to prove

\begin{theorem}
\emph{[MSh]} A profinite group is positively finitely generated if and only if
it has polynomial maximal subgroup growth.
\end{theorem}

The class of profinite groups with PMSG is very wide. For example, Borovik,
Pyber and Shalev [BPS] have shown that if the profinite group $G$ is finitely
generated, then $G$ has PMSG unless $G$ involves \emph{every finite group} as
an upper section; also iterated wreath products of finite simple groups, of
the type discussed in \S \ref{fg} above, have PMSG. So one may say that
finitely generated profinite groups have a tendency to be PFG. But the two
conditions are certainly not equivalent, since for example a non-abelian
finitely generated free profinite group (the profinite completion of a free
group) is never PFG.

Probabilistic arguments of the type given above yield all sorts of
information. The arguments always take place in the context of a profinite
group, but the conclusions sometimes apply to groups in general. I will
mention three results, all due to Mann; for the (remarkably simple) proofs,
and more discussion of the topic in general, see Chapter 11 of [SG].\medskip

\textbf{1.} Let $a_{n,d}(G)$ denote the number of $d$-generator subgroups of
index $n$ in a group $G$.

\begin{theorem}
Let $m,d\in\mathbb{N}.$ Suppose that $G$ is a group that does not involve
$\mathrm{Alt}(m+1)$ as an upper section. Then there exist $C$ and $k$,
depending only on $d$ and $m$, such that
\[
a_{n,d}(G)\leq Cn^{k}%
\]
for all $n$.
\end{theorem}

\textbf{2.} Let $\mathrm{d}(H)$ denote the minimal size of a (topological)
generating set for the profinite group $H$.

\begin{theorem}
Let $G$ be a profinite group with PSG. Then there exists a constant $C$ such
that
\[
\mathrm{d}(H)\leq C\sqrt{\log\left\vert G:H\right\vert }%
\]
for every open subgroup $H$ of $G$.
\end{theorem}

\textbf{3.} Let $h(n,r)$ denote the number of (isomorphism types of) groups of
order $n$ having a finite presentation with $r$ relations.

\begin{theorem}
Let $p$ be a prime and $r\in\mathbb{N}$. Then
\[
h(p^{k},r)=o(p^{kr})\qquad\text{as }k\rightarrow\infty.
\]

\end{theorem}

Many other results and problems are given in [M1] and [M2]. One of the most
intriguing of those was the following question: \emph{does every open subgroup
of a PFG group have PFG}? In very recent and significant work [J-ZP],
Jaikin-Zapirain and Pyber have proved that the answer is \emph{yes}. They do
this by providing a detailed characterization of groups with PMSG in terms of
the structure of their finite quotients.

\section{Other topics}

\subsection{The congruence subgroup problem\label{csp}}

I have referred to the `congruence subgroup property' in several of the
preceding sections. Recall that a subgroup $\Gamma$ in some profinite group
$G$ is said to have the CSP if the topology of $G$ induces on $\Gamma$ its own
profinite topology. This is equivalent to saying that the natural map
$\widehat{\Gamma}\rightarrow G$ is injective, or in down-to-earth terms that
\emph{every subgroup of finite index in} $\Gamma$ \emph{contains} $\Gamma\cap
N$ \emph{for some open subgroup} $N$ \emph{of} $G$. This terminology
originates in a very classical problem: what are the subgroups of finite index
in $\Gamma=\mathrm{SL}_{n}(\mathbb{Z})$? There are some obvious ones: for an
integer $m\neq0$ the \emph{principal congruence subgroup} $\operatorname{mod}%
m$ is
\begin{align*}
\Gamma(m)  &  =\left\{  g\in\Gamma\mid g_{ij}\equiv\delta_{ij}%
\,(\operatorname{mod}m)\text{ for }1\leq i,j\leq n\right\} \\
&  =\ker\left(  \Gamma\rightarrow\mathrm{SL}_{n}(\mathbb{Z}/m\mathbb{Z}%
)\right)  ,
\end{align*}
and one calls any subgroup of $\Gamma$ that contains $\Gamma(m)$ for some
$m\neq0$ a \emph{congruence subgroup}. Evidently, the congruence subgroups
have finite index in $\Gamma$, and the problem is: are there any others? This
was solved in the 1960s by Mennicke and Bass, Lazard and Serre: they proved
that the answer is `no' when $n\geq3$; as for the case $n=2$, it had been
known since the 19th century that $\mathrm{SL}_{2}(\mathbb{Z})$ has an
abundance of non-congruence subgroups of finite index.

If every subgroup of finite index is a congruence subgroup, the group $\Gamma$
is said to have the \emph{congruence subgroup property}. We see that this is a
special case of the previous definition if we consider $\Gamma$ as a subgroup
of the profinite group
\[
\widetilde{\Gamma}=\mathrm{SL}_{n}(\widehat{\mathbb{Z}}),
\]
so the congruence subgroup problem can be formulated as: is the natural map
$\widehat{\Gamma}\rightarrow\widetilde{\Gamma}$ injective?

Now $\mathrm{SL}_{n}(\mathbb{Z})$ is just the most familiar example of the
important class of $S$-\emph{arithmetic groups}, and the analogous question
applies to all such groups. I will not define these here in full generality:
for a comprehensive account see the book [PR]. Typical examples are groups of
the form $\Gamma=\mathfrak{G}(\mathbb{Z}_{S})$ where $\mathfrak{G}$ is an
algebraic matrix group defined over $\mathbb{Q}$, $S$ is a finite set of
primes and $\mathbb{Z}_{S}=\mathbb{Z}[\frac{1}{p};\,p\in S]$. The congruence
subgroup problem now becomes: determine the kernel $C(\mathfrak{G},S)$ of the
natural map
\[
\widehat{\Gamma}\rightarrow\mathfrak{G}(\widehat{\mathbb{Z}_{S}}).
\]
This group $C(\mathfrak{G},S)$ is called the \emph{congruence kernel} \ It was
observed by Serre that the natural dichotomy seems to be between those groups
whose congruence kernel is \emph{finite} and those for which it is
\emph{infinite}, and following his insight it is usual now to say that
$\Gamma$ has the CSP if $C(\mathfrak{G},S)$ is \emph{finite} (note that
according to the original definition, we would require $C(\mathfrak{G},S)=1$).
The following very general conjecture was made by Serre:\medskip

\noindent\textbf{Conjecture} Let $\mathfrak{G}$ be a simple simply connected
algebraic group over a global field $k$ and let $S$ be a finite set of places
of $k$. Then (under certain natural assumptions) the $S$-arithmetic group
$\mathfrak{G}(\mathcal{O}_{S})$ has the CSP if and only the $S$-rank of
$\mathfrak{G}$ is at least $2$.\medskip

\noindent Here, $\mathcal{O}_{S}$ denotes the ring of `$S$-integers' of $k$;
the `$S$-rank' of $\mathrm{SL}_{n}(\mathbb{Z}_{S})$, for example, is equal to
$n-1+\left|  S\right|  $. This conjecture has been proved in the majority of
cases, but some hard problems remain open: see for example [Ra].

An interesting recent development relates the congruence subgroup property of
$\Gamma$ to purely group-theoretic properties of $\Gamma$, such as its
subgroup growth and its index growth. These results are due in the main to
Platonov, Rapinchuk and Lubotzky; for a detailed account of some of them see
Chapter 7 of [SG].

\subsection{Profinite presentations\label{fp}}

By a \emph{presentation} of a group $G$ is meant an epimorphism $\pi
:F\rightarrow G$, where $F$ is a free group, together with a specific choice
of a set $X$ of free generators for $F$ and a set $R$ of generators for the
kernel $\ker\pi$ as a normal subgroup of $F$. It is usual to write
\[
G=\left\langle X\,;\,R\right\rangle ,
\]
where $R$ is a set of words on the alphabet $X$, and to interpret the symbols
in $X$ as generators of $G$ that satisfy the relations $w(X)=1$ for all $w\in
X$. For profinite groups, it is natural to consider instead epimorphisms from
a \emph{free profinite} group. When $X$ is a finite set (the only case we
consider here), the free profinite group $\widehat{F}(X)$ on $X$ is just the
profinite completion of the free group on $X$, and it has the expected
universal property with respect to continuous mappings from $X$ into profinite
groups. A \emph{profinite presentation} of $G$ is thus a \emph{continuous}
epimorphism $\pi:\widehat{F}(X)\rightarrow G$, together with a choice $R$ of
generators for $\ker\pi$ as a \emph{closed} normal subgroup of $\widehat
{F}(X)$. The elements of $R$ need no longer be words in the generators $X$: in
general they are `profinite words', that is, limits of convergent sequences of
ordinary words. But we still write
\[
G=\left\langle X\,;\,R\right\rangle
\]
to indicate such a profinite presentation (as long as the context makes it
clear which kind of presentation is meant).

The usefulness of this concept lies in the simple observation that if
$\Gamma=\left\langle X\,;\,R\right\rangle $ is an ordinary presentation of
some abstract group $\Gamma$, then $G=\left\langle X\,;\,R\right\rangle $ is a
profinite presentation of the profinite completion $G=\widehat{\Gamma}$. Given
information about a presentation of $\Gamma$, we can therefore interpret it as
information about $\widehat{\Gamma}$; profinite group theory may then yield
conclusions about $\widehat{\Gamma}$, which in turn gives us information about
$\Gamma$. This will be illustrated below. First I want to mention a celebrated
open problem.

Write $\mathrm{d}(G)$ to denote the minimal number of generators required for
a group $G$ (topological generators in the profinite context), and call
$G=\left\langle X\,;\,R\right\rangle $ a `minimal presentation' (in either
case) if $\left|  X\right|  =\mathrm{d}(G)$. The minimal number of relations
required for some minimal presentation of $G$ (in either context) is denoted
$\mathrm{t}(G)$. Now suppose that $\Gamma$ happens to be a \emph{finite}
group. In this case, of course, $\widehat{\Gamma}=\Gamma$, and we may
interpret the expression $\Gamma=\left\langle X\,;\,R\right\rangle $ either as
an ordinary presentation or as a profinite presentation. Since the topology on
$\Gamma$ is discrete, a set $X$ generates $\Gamma$ if and only if it generates
$\Gamma$ topologically. But the topology on $\widehat{F}(X)$ is by no means
discrete: just for now, let us understand $\mathrm{t}(\Gamma) $ in the
abstract sense, and write $\mathrm{t}(\widehat{\Gamma})$ for the minimal
number of relations in a minimal profinite presentation of $\Gamma$.\medskip

\noindent\textbf{Problem }Let $\Gamma$ be a finite group. Is $\mathrm{t}%
(\widehat{\Gamma})$ necesarily equal to $\mathrm{t}(\Gamma)$?\medskip

\noindent(If $r$ ordinary relations suffice to define $\Gamma,$ then the same
relations also define $\Gamma$ as a profinite group; but it is conceivable
that $\Gamma$ could be defined using a \emph{smaller} number of
\emph{profinite} relations.) For some discussion, and alternative
formulations, of this problem see \S 2.3 of [SG] (\emph{Remark} on page 48).

Two striking applications of the philosophy outlined above were made by
Lubotzky. The first uses \emph{pro-}$p$ \emph{presentations} rather than
profinite ones: these are defined in exactly the same way, using free pro-$p$
groups in place of free profinite groups.

\medskip

\textbf{1.} The famous theorem of \textbf{Golod and Shafarevich} asserts that
if $G$ is a finite $p$-group, then
\begin{equation}
\mathrm{t}(G)\geq\frac{\mathrm{d}(G)^{2}}{4} \label{gsh}%
\end{equation}
(this is true in either interpretation of the symbols, abstract or pro-$p$).
This was generalized (by Koch and Lubotzky, using Lazard's theory) to the case
of any $p$\emph{-adic analytic pro-}$p$ \emph{group }$G$ (with $\mathrm{d}(G)$
and $\mathrm{t}(G)$ now defined in terms of pro-$p$ presentations, of course).
This has consequences for any abstract group $\Gamma$ whose pro-$p$ completion
is such a group $G$; in general, $\mathrm{d}(\Gamma)$ may be strictly larger
than $\mathrm{d}(G)$, but if, for example, $\Gamma$ is nilpotent then there
exist primes $p$ such that $\mathrm{d}(\Gamma)=\mathrm{d}(\widehat{\Gamma}%
_{p})$, and one may deduce

\begin{theorem}
Let $\Gamma$ be a finitely generated non-cyclic nilpotent group. Then
$\mathrm{t}(\Gamma)\geq\mathrm{d}(\Gamma)^{2}/4.$
\end{theorem}

\noindent This is a direct generalization of the original Golod-Shafarevich
theorem to infinite groups. For details of the argument see [DDMS], Interlude
D. By further generalizing the Golod-Shafarevich theorem to a larger class of
pro-$p$ groups, J. S. Wilson established a result of still wider applicability
(it includes all finitely generated soluble groups, for example):

\begin{theorem}
\emph{[W1]} Let $\Gamma$ be a group which has no infinite $p$-torsion
residually finite quotient, for any prime $p$. Suppose that $\Gamma$ has a
presentation with $n$ generators and $r$ relations. Then%
\[
r\geq n+\frac{d^{2}-1}{4}-d
\]
where $d=\mathrm{d}(\Gamma^{\mathrm{ab}})$.
\end{theorem}

\noindent Here $\Gamma^{\mathrm{ab}}=\Gamma/\Gamma^{\prime}$ denotes the
abelianization of $\Gamma$; this appears because $\mathrm{d}(\Gamma
^{\mathrm{ab}})$ (unlike $\mathrm{d}(\Gamma)$) can be recognised as
$\mathrm{d}(\widehat{\Gamma}_{p})$ for a suitable prime $p$.

Lubotzky was concerned with groups that are very far from soluble. Let
$\Gamma$ be an \emph{arithmetic lattice} in $\mathrm{SL}_{2}(\mathbb{C})$ --
examples include groups like $\mathrm{SL}_{2}(\mathcal{O})$ where
$\mathcal{O}$ is the ring of integers in an imaginary quadratic field, but
there are more mysterious ones. It is fairly easy to see that if $\Gamma$ has
the congruence subgroup property then its pro-$p$ completion $\widehat{\Gamma
}_{p}=G$ is $p$-adic analytic, and hence satisfies (\ref{gsh}); moreover, the
same holds for the pro-$p$ completion of every subgroup $\Delta$ of finite
index in $\Gamma$. From this it may be deduced that
\[
\left|  X\right|  -\left|  R\right|  \leq\mathrm{d}_{p}(\Delta)-\frac
{\mathrm{d}_{p}(\Delta)^{2}}{4}%
\]
for every finite presentation $\Delta=\left\langle X\,;\,R\right\rangle $,
where
\[
\mathrm{d}_{p}(\Delta)=\mathrm{d}(\widehat{\Delta}_{p})=\mathrm{d}%
(\Delta/\Delta^{p}[\Delta,\Delta])\text{.}%
\]
On the other hand, according to a theorem of Epstein each such $\Delta$ has a
presentation $\left\langle X\,;\,R\right\rangle $ for which $\left|  R\right|
\leq\left|  X\right|  $ (assuming, as we may, that $\Delta$ is torsion-free).
Hence $\mathrm{d}_{p}(\Delta)\leq4$. Now the theory of linear groups shows
that if the numbers $\mathrm{d}_{2}(\Delta)$ are \emph{bounded} as $\Delta$
ranges over all the subgroups of finite index in $\Gamma$, then $\Gamma$ must
have a soluble subgroup of finite index. This is certainly not the case here,
so we have

\begin{theorem}
\emph{[Lu1]} No arithmetic lattice in $\mathrm{SL}_{2}(\mathbb{C})$ satisfies
the congruence subgroup property.
\end{theorem}

\noindent This establishes many of the `negative' cases of Serre's conjecture,
stated in the preceding subsection. The method has been generalized by
Lubotzky to obtain

\begin{theorem}
Let $\Gamma$ be any lattice in $\mathrm{SL}_{2}(\mathbb{C})$. Then $\Gamma$
has subgroup growth of type at least $n^{(\log n)^{2-\varepsilon}}$ for every
$\varepsilon>0$.
\end{theorem}

A \emph{lattice} is a discrete subgroup of finite co-volume. Since the
\emph{congruence subgroup growth} of any arithmetic group is at most of type
$n^{\log n/\log\log n}$, this shows that the congruence subgroup property
fails here in a dramatic way: the subgroups of finite index vastly outnumber
the congruence subgroups as the index goes to infinity. For details of the
proof, and many other cases, see Chapter 7 of [SG].

Returning to profinite groups, or rather pro-$p$ groups, the most powerful
generalization of the Golod-Shafarevich theorem was obtained by Zelmanov:

\begin{theorem}
\emph{[Z]} Let $G$ be a non-procyclic finitely generated pro-$p$ group with a
minimal pro-$p$ presentation $G=\left\langle X\,;\,R\right\rangle $. Then
either $\left|  R\right|  \geq\left|  X\right|  ^{2}/4$ or else $G$ contains a
closed subgroup that is a non-abelian free pro-$p$ group.
\end{theorem}

\textbf{2.} Let $f(n,d)$ denote the number of (isomorphism types of)
$d$-generator groups of order $n$. Establishing a conjecture of Mann, Lubotzky proved

\begin{theorem}
For every $n$ and $d$ we have
\[
f(n,d)\leq n^{2(d+1)\lambda(n)}.
\]

\end{theorem}

\noindent Here $\lambda(n)=\sum l_{i}$ where $n=\prod p_{i}^{l_{i}}$ is the
factorization of $n$ into prime-powers. This is deduced from the following
theorem: \emph{every finite simple group of order} $n$ \emph{has a profinite
presentation with }$2$ \emph{generators and at most} $2\lambda(n)$
\emph{relations}. It is conjectured that this remains true if the word
`profinite' is omitted, and this has been proved in most cases. But it is in
general easier to get at a profinite presentation than at an abstract
presentation: roughly speaking, if $N=\ker\pi$ in our original notation, then
the number of \emph{profinite} relations needed for a presentation
$\pi:F\rightarrow G$ can be detected in the `relation module' $N/[N,N]$,
whereas the number of `ordinary' relations depends on the structure of $N$
itself as an $F$-operator group. For details, see \S 2.3 of [SG].

\subsection{Profinite trees\label{free}}

A large part of combinatorial group theory deals with the properties of
generalized free products and HNN extensions. A powerful unified framework for
studying such constructions is the Bass-Serre theory of groups acting on
trees. In recent years, an analogous theory has been developed of profinite
groups acting on `profinite trees', largely due to the work of Melnikov, Ribes
and P. A. Zalesskii. As well as providing a basis for the theory of
generalized free products in the profinite category, this has found a number
of applications to to abstract free groups and free products; a typical
example is Theorem \ref{csep} mentioned in \S \ref{local}, above.

This is a significant chapter in `pure' profinite group theory, with solid
achievements but also presenting a number of challenging open problems.
However, it is beyond my competence to present anything like an adequate
account of it. Detailed expositions of the theory are given in [RZ1] (for
pro-$p$ groups) and the forthcoming book [RZ3]; for various specific
applications, see the papers [RZ4], [RZ5], [RZ6] and [RSZ].

\bigskip

\begin{center}
{\Large References}

\bigskip
\end{center}

[BG] L. Bartholdi, R. I. Grigorchuk and Z. \u{S}uni\'{k}, Branch groups, in
\emph{Handbook of Algebra III}, ed. M.Hazewinkel, North-Holland, Amsterdam, 2003.

\bigskip

[B] N. Boston, $p$-adic Galois representations and pro-$p$ Galois groups.
\emph{Chapter 11 in }[NH].

\bigskip

[BCRS] G. Baumslag, F. Cannonito, D. J. S. Robinson and D. Segal, The
algorithmic theory of polycyclic-by-finite groups, \emph{J. Algebra}
\textbf{142} (1991), 118-149.

\bigskip

[BMP] A. Balog, A. Mann and L. Pyber, Polynomial index growth groups,
\emph{Int. J. Algebra and Computation} \textbf{10} (2000), 773-782.

\bigskip

[BPS] A. Borovik, L. Pyber \& A. Shalev, Maximal subgroups in finite and
profinite groups, \emph{Trans. Amer. Math. Soc}. \textbf{348} (1996), 3745-3761.

\bigskip

[CG] J-P. Serre, \emph{Galois Cohomology.} Springer Verlag, Berlin-Heidelberg, 1997.

\bigskip

[DDMS] J. D. Dixon, M. P. F. du Sautoy, A. Mann \& D. Segal, \emph{Analytic
pro-}$p$ \emph{Groups}, 2nd edition, Cambridge Studies in Advanced Maths.
\textbf{61}, Cambridge Univ. Press, Cambridge, 1999.

\bigskip

[dS1] M. P. F. du Sautoy, Finitely generated groups, $p$-adic analytic groups
and Poincar\'{e} series, \emph{Annals of Math} \textbf{137} (1993), 639-670.

\bigskip

[dS2] M. P. F. du Sautoy, Counting $p$-groups and nilpotent groups,
\emph{Publ. Math. IHES} \textbf{92} (2000), 63-112.

\bigskip

[dSF] M. P. F. du Sautoy and I. Fesenko, Where the wild things are:
ramification groups and the Nottingham group. \emph{Chapter 10 in
}[NH]\emph{.}

\bigskip

[dSS] M. P. F. du Sautoy and D. Segal, Zeta functions of groups. \emph{Chapter
9 in }[NH].

\bigskip

[E1] B. Eick, Orbit-stabilizer problems and computing normalizers for
polycyclic groups, \emph{J. Symbolic Comp}. \textbf{34} (2002), 1-19.

\bigskip

[E2] B. Eick, Computing with infinite polycyclic groups. \emph{Groups and
computation, III (Columbus, Ohio 1999)}, pp. 139-154, Ohio State Univ. Math.
Res. Inst. Publ. 8, de Gruyter, Berlin 2001.

\bigskip

[FJ] M. D. Fried \& M. Jarden, \emph{Field arithmetic}, Ergebnisse der Math.
(3) \textbf{11}, Springer-Verlag, Berlin--Heidelberg, 1986.

\bigskip

[G] R. I. Grigorchuk, Just infinite branch groups. Chapter 4 in [NH].

\bigskip

[GPS] F. J. Grunewald, P. F. Pickel and D. Segal, Polycyclic groups with
isomorphic finite quotients, \emph{Annals of Math. }\textbf{111} (1980), 155-195.

\bigskip

[GSS] F. J. Grunewald, D. Segal \& G. C. Smith, Subgroups of finite index in
nilpotent groups,\emph{ Invent. Math.} \textbf{93} (1988), 185-223.

\bigskip

[H] D. F. Holt, Enumerating perfect groups, \emph{J. London Math. Soc}.
\textbf{39} (1989), 67-78.

\bigskip

[J-Z1] A. Jaikin-Zapirain, Zeta function of representations of compact
$p$-adic analytic groups, \emph{J. Amer. Math. Soc.} \textbf{19} (2006), 91-118.

\bigskip

[J-Z2] A. Jaikin-Zapirain, On the verbal width of finitely generated pro-$p$
groups, \emph{in preparation.}

\bigskip

[J-ZP] A. Jaikin-Zapirain and L. Pyber, Random generation of finite and profinite groups and group enumeration,
\emph{to appear.}

\bigskip

[KN] M. Kassabov and N. Nikolov, Cartesian products as profinite completions,
\emph{International Math. Research Notices}  \textbf{2006} (2006), Article ID 72947.

\bigskip

[L] M. Lazard, Groupes analytiques $p$-adiques, \emph{Publ. Math. IHES}
\textbf{26} (1965), 389-603.

\bigskip

[LL] M. Larsen and A. Lubotzky, Representation Growth for Linear Groups,
\texttt{arXiv:math.GR/0607369}.

\bigskip

[LGM1] C. R. Leedham-Green and S. McKay, \emph{The structure of groups of
prime-power order}, LMS Monographs \textbf{27}, Oxford Univ. Press, Oxford, 2002.

\bigskip

[LGM2] C. R. Leedham-Green and S. McKay, On the classification of $p$-groups
and pro-$p$ groups. \emph{Chapter 2 in }[NH]\emph{.}

\bigskip

[Lu1] A. Lubotzky, Group presentation, $p$-adic analytic groups and lattices
in $\mathrm{SL}_{2}(\mathbb{C})$. \emph{Ann. of Math.} (2) \textbf{118}
(1983), 115--130.

\bigskip

[Lu2] A. Lubotzky, A group-theoretic characterization of linear groups,
\emph{J. Algebra }\textbf{113} (1988), 207-214.

\bigskip

[LM] A. Lubotzky \& A. Mann, Powerful $p$-groups. I: Finite groups; II:
$p$-adic analytic groups, \emph{J. Algebra} \textbf{105} (1987), 484-505, 506-515.

\bigskip

[LMS] A. Lubotzky, A. Mann \& D. Segal, Finitely generated groups of
polynomial subgroup growth, \emph{Israel J. Math.} \textbf{82} (1993), 363-371.

\bigskip

[LS] R. C. Lyndon and P. E. Schupp, \emph{Combinatorial group theory,}
Springer-Verlag, Berlin-Heidelberg-New York, 1977.

\bigskip

[LW] P. A. Linnell and D. Warhurst, Bounding the number of generators of a
polycyclic group, \emph{Archiv der Math.} \textbf{37} (1981), 7-17.

\bigskip

[M1] A. Mann, Positively finitely generated groups, \emph{Forum Math.}
\textbf{8} (1996), 429-459.

\bigskip

[M2] A. Mann, Some applications of probability in group theory, \emph{Groups:
topological, combinatorial and arithmetic aspects}, ed. T. W. M\"{u}ller, LMS
Lect. note series \textbf{311}, CUP, Cambridge, 2004.

\bigskip

[MS1] A. Mann and D. Segal, Uniform finiteness conditions in residually finite
groups, \emph{Proc. London Math. Soc.} (3) \textbf{61} (1990), 529-545.

\bigskip

[MS2] A. Mann and D. Segal, Subgroup growth: some current developments,
in\emph{ Infinite Groups 94}, eds. de Giovanni and Newell, W. de Gruyter, 1995.

\bigskip

[MSh] A. Mann and A. Shalev, Simple groups, maximal subgroups and
probabilistic aspects of profinite groups, \emph{Israel J. Math. }\textbf{96}
(1996), 449-468.

\bigskip

[NH] M. P. F. du Sautoy, D. Segal \& A. Shalev (ed.), \emph{New horizons in
pro-}$p$\emph{ groups,} Progress in Math. \textbf{184}, Birkh\"{a}user Boston, 2000.

\bigskip

[NS1] N. Nikolov and D. Segal, Finite index subgroups in profinite groups,
\emph{C. R. Acad. Sci. Paris, Ser. I }\textbf{337} (2003), 303-308.

\bigskip

[NS2] N. Nikolov and D. Segal, On finitely generated profinite groups, I:
strong completeness and uniform bounds. II: product decompositions of
quasisimple groups, \emph{Annals of Math. }\textbf{165} (2007), 171-238, 239-273.

\bigskip

[NS3] N. Nikolov and D. Segal, Direct products and profinite completions,
\emph{J. Group Theory,} to appear.

\bigskip

[PR] V. P. Platonov \& A. S. Rapinchuk, \emph{Algebraic groups and number
theory}, Academic Press, San Diego, 1994.

\bigskip

[P] L. Pyber, Groups of intermediate subgroup growth and a problem of
Grothendieck, \emph{Duke Math. J. }\textbf{121 }(2004), 169-188.

\bigskip

[PS] L. Pyber and D. Segal, Finitely generated groups with polynomial index
growth, \emph{J. reine angewandte Math., }to appear.

\bigskip

[Ra] A. S. Rapinchuk, The congruence subgroup problem. \emph{Algebra,
K-theory, groups, and education} (New York, 1997), 175--188, \emph{Contemp.
Math.} \textbf{243}, Amer. Math. Soc., Providence, RI, 1999.

\bigskip

[R] V. A. Romankov, Width of verbal subgroups in solvable groups,
\emph{Algebra i Logika} \textbf{21}(1982), 60-72 (Russian);\textbf{ }
\emph{Algebra and Logic} \textbf{21 }(1982), 41-49 (English).

\bigskip

[RZ1] L. Ribes and P. A. Zalesskii, Pro-$p$ trees and applications.
\emph{Chapter 3 in }[NH].

\bigskip

[RZ2] L. Ribes and P. A. Zalesskii, \emph{Profinite groups}. Ergebnisse der
Math. \textbf{40}, Springer, Berlin -- Heidelberg , 2000.

\bigskip

[RZ3] Ribes and Zalesskii, \emph{Profinite trees}, Springer, Berlin (to appear).

\bigskip

[RZ4] Ribes and Zalesskii, On the profinite topology on a free group,
\emph{Bull. London Math. Soc.} \textbf{25} (1993), 37-43.

\bigskip

[RZ5] Ribes and Zalesskii, The pro-$p$ topology of a free group and
algorithmic problems in semigroups, \emph{Int. J. Algebra and Comput.}
\textbf{4} (1994), 359-374.

\bigskip

[RZ6] Ribes and Zalesskii, Profinite topologies in free products of groups,
\emph{Int. J. Algebra and Comput}. \textbf{14} (2004), 751-772.

\bigskip

[RSZ] L. Ribes, D. Segal and P. A. Zalesskii, Conjugacy separability and free
products with cyclic amalgamation, \emph{J. London Math. Soc.} (2) \textbf{57}
(1998), 609-628.

\bigskip

[S] D. Segal, \emph{Polycyclic groups}, Cambridge Univ. Press, Cambridge,
1983. (\emph{Reprinted in paperback, }2005.)

\bigskip

[S1] D. Segal, Decidable properties of polycyclic groups, \emph{Proc. London
Math. Soc.} (3) \textbf{61} (1990), 497-528.

\bigskip

[S2] D. Segal, Variations on polynomial subgroup growth, \emph{Israel J.
Math.} \textbf{94} (1996), 7-19.

\bigskip

[S3] D. Segal, Subgroups of finite index in soluble groups II, in \emph{Groups
St Andrews 1985}, LMS Lect. note series \textbf{121}, pp. 315-319, CUP,
Cambridge, 1986.

\bigskip

[S4] D. Segal, The finite images of finitely generated groups, \emph{Proc.
London Math. Soc}. (3) \textbf{82} (2001), 597-613.

\bigskip

[Sh] A. Shalev, The structure of finite $p$-groups: effective proof of the
coclass conjectures, \emph{Invent. Math.} \textbf{115} (1994), 315-345.

\bigskip

[SG] A. Lubotzky and D. Segal, \emph{Subgoup Growth,} Progress in Math.
\textbf{212}, Birkh\"{a}user, Basel, 2003.

\bigskip

[SSh] D. Segal \& A. Shalev, Profinite groups with polynomial subgroup growth,
\emph{J. London Math. Soc}.(2) \textbf{55} (1997), 320-334.

\bigskip

[V] C. Voll, Functional equations for zeta functions of groups and rings,
\emph{preprint}.

\bigskip

[W1] J. S. Wilson, Finite presentations of pro-$p$ groups and discrete groups,
\emph{Invent. Math.} \textbf{105} (1991), 177-183.

\bigskip

[W2] J. S. Wilson, \emph{Profinite groups}, London Math. Soc. Monographs
(n.s.) \textbf{19}, Clarendon Press, Oxford, 1998.

\bigskip

[Z] E. Zelmanov, On groups satisfying the Golod-Shafarevich condition.
\emph{Chapter 7 in }[NH].

\bigskip

\texttt{All Souls College}

\texttt{Oxford OX1 4AL}

\texttt{UK.}
\end{document}